\documentstyle[11pt,ctagsplt]{amsart}

\newtheorem{thm}{Theorem}
\newtheorem{prop}{Proposition}
\newtheorem{lemma}{Lemma}
\newtheorem{defn}{Definition}

\newcommand{\R}{\Bbb R}

\newcommand{\bs}{\backslash}
\newcommand{\norm}[1]{\left|#1\right|}

\newcommand{\HR}{\Bbb H_{\Bbb R}^2}
\newcommand{\HC}{\Bbb H_{\Bbb C}^2}
\newcommand{\Div}{\operatorname{div}}

\newcommand{\M}{(M,g)}
\newcommand{\vp}{\varphi}
\newcommand{\D}{\partial}
\newcommand{\G}{\Gamma}
\newcommand{\dist}{\operatorname{dist}}
\newcommand{\compose}{\cdot}
\renewcommand{\S}{\cal S}
\newcommand{\B}{\cal B}
\newcommand{\disp}{\displaystyle}
\renewcommand{\H}{\cal H}
\newcommand{\Hc}{\cal H_R}

\newcommand{\xib}{\overline\xi}
\newcommand{\phib}{\overline\phi}
\newcommand{\ub}{\overline u}
\newcommand{\vb}{\overline v}
\newcommand{\Qb}{\overline Q}
\newcommand{\Tb}{\overline T}
\newcommand{\chib}{{\overline\chi}}
\newcommand{\sphere}{\Bbb S}
\newcommand{\Ha}{\Bbb H}
\newcommand{\Co}{\Bbb C}
\newcommand{\K}{\Bbb K}
\newcommand{\inner}[2]{\langle#1,#2\rangle}
\newcommand{\bx}{\bold x}
\newcommand{\xb}{\overline x}
\newcommand{\by}{\bold y}
\newcommand{\bz}{\bold z}
\newcommand{\bw}{\bold w}
\newcommand{\scalar}[2]{(#1,#2)}
\newcommand{\zb}{\overline z}
\renewcommand{\o}{{\bold o}}
\newcommand{\Do}{\Bbb D}
\renewcommand{\P}{\Bbb P}
\newcommand{\e}{\bold e}
\renewcommand{\Re}{\operatorname{Re}}
\renewcommand{\Im}{\operatorname{Im}}

\newcommand{\aut}[1]{{\sc #1}}
\newcommand{\tit}[1]{{\em #1\/}}
\newcommand{\vol}[1]{{\bf #1}}
\newcommand{\yr}[1]{(#1)}
\newcommand{\pp}[2]{#1--#2}

\begin{document}

\title[Harmonic Maps with Prescribed Singularities]
	{On the Dirichlet Problem for\\
	Harmonic Maps with Prescribed Singularities}
\author{Gilbert Weinstein}
\address{University of Alabama at Birmingham, Birmingham, Alabama 35294}
\email{weinstei@@math.uab.edu}
\thanks{Research supported in part by
	a grant from NSF EPSCoR in Alabama}
\date{\today}

\maketitle

\begin{abstract}
Let $\M$ be a classical Riemannian globally symmetric space
of rank one and non-compact type.  We prove the existence and
uniqueness of solutions to the Dirichlet problem for harmonic maps into
$\M$ with prescribed singularities along a closed submanifold of the
domain.
This generalizes our previous work where such maps into the
hyperbolic plane were constructed.
This problem, in the case where
$\M$ is the complex-hyperbolic plane, has
applications to equilibrium configurations of co-axially
rotating charged black holes in General Relativity.
\end{abstract}

\section{Introduction}

The Einstein vacuum equations in the stationary axially symmetric
case reduce to a harmonic map from $\R^3$
into $\HR$, the hyperbolic plane, with prescribed
singularities along the axis of symmetry.
In~\cite{weinstein90,weinstein92}, we used this fact to construct solutions
of these equations which could be interpreted as a pair of rotating
black holes held apart by a singular strut.  These solutions generalized
the static Weyl solutions, see~\cite{weyl}.
The first step in this program was to solve a Dirichlet problem for such
maps with the singularity prescribed along a closed submanifold of the
domain.  A
natural generalization of this problem is to replace the Einstein vacuum
equations with the Einstein-Maxwell equations.
A similar reduction again leads to a harmonic map problem with prescribed
singularities, but the target is now $\HC$,
the complex hyperbolic plane, see~\cite{mazur}.

In this paper, we study the Dirichlet problem for harmonic maps with
prescribed singularities from a smooth bounded domain $\Omega\subset\R^n$,
$n\geq2$, into $\M$ a classical Riemannian globally
symmetric space of rank one and of non-compact type.  Thus $\M$ is either
the real-, complex-, or quaternion-hyperbolic space,
i.e.\ $\M = \Ha^\ell_{\K}$, where $\ell\geq2$, and
$\K$ is either $\R$, $\Co$, or the quaternions $\Ha$,
see~\cite{helgason}.
For simplicity, we take the Euclidean metric on $\R^n$, although all the
results carry over easily to bounded domains in Riemannian manifolds.
Recall that a map $\vp\colon\Omega\to\M$ is {\em harmonic\/} if for each
$\Omega'\subset\subset\Omega$ the map $\vp|\Omega'$ is
a critical point of the energy:
\begin{equation}
	\label{energy}
	E_{\Omega'}(\vp) = \int_{\Omega'} \norm{d\vp}^2,
\end{equation}
where $\norm{d\vp}^2 =
\sum_{k=1}^n g(\nabla_k\vp,\nabla_k\vp)$.  It then satisfies an
elliptic system of nonlinear partial differential equations, written in
local coordinates on $M$ as:
\[
	\Delta \vp^a +
	\sum_{k=1}^n \G^a_{bc} \, \D_k\vp^b \, \D_k\vp^c = 0,
\]
where $\G^a_{bc}$ are the Christoffel symbols of $\M$.  Harmonic maps have
been studied extensively.  The Dirichlet problem for harmonic maps into a
manifold of non-positive curvature was first solved by R.~Hamilton
in~\cite{hamilton} using a
heat flow method.  A variational approach was later developed by
R.~Schoen and K.~Uhlenbeck, see~\cite{schoen82,schoen83}.
More recently, P.~Li and L.-F.~Tam constructed harmonic maps between
hyperbolic spaces, see~\cite{lt91,lt93}.

It is well known that if $\M$ has negative sectional curvature and
$\vp\colon\Omega\to\M$ is a finite energy harmonic map
then $\vp\in C^\infty(\Omega;M)$.
Furthermore, if $\D\Omega$ is of class $C^{2,\alpha}$, and $\vp|\D\Omega$ is
$C^{2,\alpha}$, then
$\vp\in C^{2,\alpha}(\overline\Omega;M)$, see~\cite{schoen83}.

Let $\Sigma_i$, $i=1,\dots N$, be disjoint
closed smooth submanifolds of $\Omega$ of co-dimension at least $2$, and
set $\Sigma=\cup_{i=1}^N \Sigma_i$.  For each $1\leq i\leq N$, let
$\gamma_i\colon\R\to\M$ be a unit speed geodesic, and let $\vp_i\colon
\Omega\bs\Sigma_i\to\M$ be a harmonic map singular on
$\Sigma_i$ whose image is contained in $\gamma_i(\R)$ and
such that $\vp(x)\to\gamma_i(+\infty)\in\D M$ as $x\to\Sigma_i$.
We shall call such a map a {\em $\Sigma_i$-singular map into $\gamma_i$},
provided it satisfies an additional technical condition.
Since $\gamma_i$ is (trivially) flat and totally geodesic,
such a map is easily constructed from a harmonic function
$u_i$ on $\Omega\bs\Sigma_i$ which tends to infinity on $\Sigma_i$.
Let $\psi\colon\D\Omega\to M$ be a smooth boundary map.
We wish to find a harmonic map $\vp\colon\Omega\bs\Sigma\to\M$ which has
boundary values $\psi$, and is asymptotic to $\vp_i$ near $\Gamma_i$, see
section~\ref{prelim} for the definitions.
The main result of this paper is the following theorem:

\begin{thm}	\label{main}
For $1\leq i\leq N$, let $\vp_i$ be a $\Sigma_i$-singular map
into $\gamma_i$, and let $\psi\in C^{2,\alpha}(\D\Omega;M)$.
Then, there exists a unique
harmonic map $\vp\in C^\infty(\Omega\bs\Sigma;M)\cap
C^{2,\alpha}(\overline\Omega\bs\Sigma;M)$ , such that
$\vp=\psi$ on $\D\Omega$, and $\vp$ is asymptotic to $\vp_i$ near
$\Sigma_i$ for each $1\leq i \leq N$.
\end{thm}

On the one hand, these maps may be viewed as non-linear generalizations of
harmonic functions $u$ which tend to $\pm\infty$ on $\Sigma$,
the case $m=1$.  Note that in this case there are only two points at
infinity in $\M$.
On the other hand, they may be
viewed as generalizations of geodesic rays, the case $n=1$.  However, in
this case $\Sigma$ is necessarily
empty, and thus singular asymptotic behavior can
only be prescribed at infinity.  Also, we should point out that M.~Anderson
constructed in~\cite{anderson}
complete area-minimizing hypersurfaces in hyperbolic
spaces asymptotic to a given set at infinity.  However, he assumes
the given set cuts the boundary of $\M$ into exactly two
connected components, a situation entirely different from ours.

We use a direct variational method to prove Theorem~\ref{main}, following
the same outline as in~\cite{weinstein90,weinstein92}.
The difficulty is that the prescribed singularities force all
admissible maps to have infinite energy on $\Omega$.
To remedy this, we renormalize the energy,
making use of the Busemann functions on $\M$.
Busemann functions were used
in~\cite{shen} to prove a Liouville type theorem for harmonic maps into
negatively curved manifolds.
Our method would apply to construct harmonic
maps with prescribed singularities into any simply connected manifold of
pinched negative curvature were it not for the fact that we use the
specific structure of $\M$ in a crucial way in Lemma~\ref{lemma:Qequiv}.
In another direction, one could try to relax the curvature condition to
allow for symmetric spaces of rank $\geq 2$ as targets.

The plan of the paper is as follows.  In Section~\ref{prelim}, we
discuss some preliminaries.
This includes a detailed study of manifolds with pinched negative curvature
necessary for the
variational approach to go through.  Also included in this section
are some definitions, and a maximum principle needed for the uniqueness.
In Section~\ref{single}, we present the proof of Theorem~\ref{main} in the
somewhat simpler case $N=1$ where only one
singular asymptotic behavior is prescribed for $\vp$.
Nevertheless, the proof of this case already contains most of the
main ideas.
Then, in Section~\ref{multiple}, we treat the
case $N\geq2$, where multiple behaviors are prescribed.
In an appendix, we provide some of the calculations needed in the proof of
Lemma~\ref{lemma:Qequiv}.

In a forthcoming paper, we will
treat the case of unbounded domains, where singular asymptotic behavior
should also be prescribed at infinity, and we will apply these results
to the rotating charged black hole problem in General Relativity,
as mentioned in the first paragraph of this introduction.

\section{Preliminaries} \label{prelim}

Let $\M$ be an $m$-dimensional
simply connected Riemannian manifold with sectional curvatures
bounded between two negative constants: $-b^2\leq\kappa\leq -a^2<0$.  Thus,
$\M$ is a Cartan-Hadamard manifold.  We first recall
a few standard facts
about this class of manifolds taken mostly from~\cite{eberlein}
and~\cite{heintze}.
Throughout, all geodesics are unit speed.
Let $\gamma\colon\R\to\M$ be a geodesic.
The Busemann function associated with $\gamma$ is defined by:
\[
	f_{\gamma}(p) = \lim_{t\to\infty}
	\bigl( \dist\bigl(p,\gamma(t)\bigr) - t\bigr).
\]
This is the {\em renormalized\/} distance function from the ideal point
$\gamma(+\infty)\in\D M$.  As such, it inherits many of the
properties of $\dist(\cdot,q)$, the distance function from a fixed point
$q$.  In particular,
it is convex, and its gradient has length $1$:
\begin{align}
	\label{norm1}
	\norm{\nabla f_{\gamma}} &= 1,\\
	\label{convex}
	\qquad \nabla^2 f_{\gamma} &\geq 0.
\end{align}
Furthermore, $f_{\gamma}\in C^2(M)$, and
the level sets $\S_\gamma(t)= \{ p\in M: f_{\gamma}(p)= t \}$,
called {\em horospheres\/}, are $C^2$-diffeomorphic to $\R^{m-1}$.

Denote the reverse geodesic $t\mapsto\gamma(-t)$ by $-\gamma$.  To be
consistent with the notation used later, we now use $-\gamma$ in place of
$\gamma$.  Let $v_0\colon\S_{-\gamma}(0)\to\R^{m-1}$
be a $C^2$-global coordinate system on $\S_{-\gamma}(0)$ centered at
$\gamma(0)$.  From~\eqref{norm1}, it follows that the integral
curves of $\nabla f_{-\gamma}$,
the field of unit normals to the horospheres, are geodesics.
Let $\phi_t$ be the flow generated by this vector field,
then $\phi_{-t}$ maps $\S_{-\gamma}(t)$ to $\S_{-\gamma}(0)$, and
$v_t=v_0\compose\phi_{-t}$ is a $C^1$-coordinate system on $\S_{-\gamma}(t)$.
Define $v\colon M\to\R^{m-1}$ by $v|\S_{-\gamma}(t) = v_t$, and let
$u=f_{-\gamma}$, then $\phi=(u,v)\colon M\to\R^m$
is a $C^1$-coordinate system on $M$.
In this coordinate system, the metric $g$ can be written as:
\[
	ds^2 =  du^2 + Q_p(dv),
\]
where, for each
$p\in M$, $Q_p$ is a positive quadratic form on $\R^{m-1}$, and $Q_p$ is
continuous in $p$.
Specific examples of this construction are given in
Lemma~\ref{lemma:Qequiv}.
{}From~\eqref{convex}, it is easily seen that for each non-zero
$\xi\in\R^{m-1}$, and each
fixed $v\in\R^{m-1}$, $Q_{\phi^{-1}(u,v)}(\xi)$ is a
positive continuous increasing function of $u$.

We now wish to sharpen this result in the lemma below.
Following~\cite{heintze},
we will say that a Jacobi field $Y$ along a geodesic $\gamma$
is {\em stable as $t\to\pm\infty$\/} if it is bounded for $\pm t\geq0$.
For each $v\in\R^{m-1}$,
let $\gamma_v$ denote the geodesic $t\mapsto\phi^{-1}(t,v)$, and note that
$\gamma_0=\gamma$.

\begin{lemma} \label{lemma:expansion}
For every $v\in\R^{m-1}$, and every $t\in\R$, there holds:
\begin{equation}	\label{expansion}
	2a \, Q_{\gamma_v(t)}(\xi) \leq
	\frac{d}{dt} \left(Q_{\gamma_v(t)}(\xi)\right)
	\leq  2b \, Q_{\gamma_v(t)}(\xi), \qquad \forall\xi\in\R^{m-1}.
\end{equation}
\end{lemma}
\begin{pf}
We first make the following observation.  Fix $\xi\in\R^{m-1}$, let
$\xib = (0,\xi)\in\R^m$, and define $Y_{\xi}=d\phi^{-1} \cdot \xib$.
Then, $Y_\xi\bot\dot\gamma_v$ everywhere, and
$Y_{\xi}$ is a Jacobi field along $\gamma_v$ which is stable
as $t\to-\infty$.  To see this it suffices to note that for each $s\in\R$,
the curve $t\mapsto\phi^{-1}(t,v+s\xi)$ is the geodesic $\gamma_{v+s\xi}$,
$Y_\xi$ is the variation vector field of this family, and
$\dist(\gamma_{v+s\xi}(t),\gamma_v(t))$ is bounded for $t\leq0$.  Since
\[
	Q_p(\xi) = g(d\phi^{-1}|_{\phi(p)} \cdot\xib,
		d\phi^{-1}|_{\phi(p)} \cdot\xib)
		= g\bigl(Y_\xi(p),Y_\xi(p)\bigr),
\]
we see that~\eqref{expansion} is simply an estimate on the
logarithmic growth rate of stable Jacobi fields in $\M$.
Fix $v\in\R^{m-1}$, and let
\[
	\chi(t) = Q_{\gamma_v(t)}(\xi) = g(Y_\xi,Y_\xi)|_{\gamma_v(t)}.
\]
In~\cite{heintze},
a pointwise estimate is proved which in our notation reads:
\begin{equation}	\label{jacobi}
	e^{2a(t-s)}\chi(s) \leq \chi(t) \leq e^{2b(t-s)}\chi(s),
\end{equation}
for all $s\leq t$.  The lemma follows
by taking logarithms, dividing by $t-s$, and letting $t\to s^+$.
\end{pf}

Two geodesics $\gamma$ and $\gamma'$
are said to be {\em asymptotic\/}
if $\dist\bigl(\gamma(t),\gamma'(t)\bigr)$
is bounded for $t\geq0$.  This is clearly an equivalence
relation.  The boundary $\D M$ of $\M$
is defined to be the set of equivalence classes of geodesics in $\M$.
We denote the equivalence class of $\gamma$ in $\D M$ by $\gamma(+\infty)$.
Since the sectional curvatures of $\M$ are pinched between
two negative constants, there is for each pair $\omega\ne\omega'\in\D M$
a geodesic from $\omega$ to $\omega'$ unique up to translation.
We also write $\gamma(-\infty)$ for $-\gamma(+\infty)$.
As a corollary of Lemma~\ref{lemma:expansion} we note the following:

\begin{lemma}	\label{asymptotic}
Let $\gamma$ and $\beta$ be geodesics in $\M$ such that
$\gamma(-\infty)=\beta(-\infty)$, parameterized so that
$f_{-\gamma}=f_{-\beta}$.  Then
\begin{equation}
	\lim_{t\to-\infty} \dist\bigl(\gamma(t),\beta(t)\bigr)=0.
\end{equation}
\end{lemma}
\begin{pf}
With the coordinate system $\phi$ defined above,
we may assume that $\beta=\gamma_v$ for some $v\in\R^{m-1}$.  For each
$t\in\R$, define the curve $\sigma_t\colon[0,1]\to M$ by $\sigma_t(s) =
\phi^{-1}(t,sv)$ for $0\leq s\leq1$.  We have $\sigma_t(0)=\gamma(t)$ and
$\sigma_t(1)=\beta(t)$, hence if $d(t)=\dist(\gamma(t),\beta(t))$,
clearly $d(t)\leq$ the length of $\sigma_t$, and therefore
\[
	\bigl(d(t)\bigr)^2 \leq \int_0^1 g(\dot\sigma,\dot\sigma) =
	\int_0^1 Q_{\sigma_t(s)} (v).
\]
Denoting the integral on the right-hand side by $h(t)$, then
Lemma~\ref{lemma:expansion} implies that
\[
	\frac{dh}{dt} \geq 2ah,
\]
which gives $h(t) \leq e^{2at} h(0)$ for all $t\leq0$.  The lemma follows.
\end{pf}

The {\em horoball\/} associated with $\gamma$ is defined
by:
\[
	\B_\gamma(t) = \{ p\in M: f_\gamma(p) \leq t\}.
\]
Denote the closed geodesic ball of radius $R>0$ centered at $p\in M$
by $B_R(p)$.

We have the following lemmas:

\begin{lemma}	\label{bound}
Let $\gamma$ be a geodesic in $\M$.  Then for any $t_0\in\R$, and any
$T\geq0$, we have
\[
	\B_{\gamma}(-t_0+T) \cap \B_{-\gamma}(t_0+T)
	\subset B_R(\gamma(t_0)),
\]
where $R= T + a^{-1} \log2$.
\end{lemma}
\begin{pf}
We first remark that, by shifting the parameter along $\gamma$, we may
without loss of generality assume that $t_0=0$.  We use the
following comparison principle which is proved in~\cite[Lemma~4.2]{heintze}:
\begin{quote}
	{\em Let $\beta$ be any geodesic in $\M$,
	and let $\chi=f_\gamma\compose\beta$.
	Similarly, let $\chi_a$ be the restriction of a Busemann function
	along a geodesic in a space of constant curvature $-a^2$.
	Suppose that $\chi(0)=\chi_a(0)$, and
	that $\dot\chi(0)=\dot\chi_a(0)$.  Then, $\chi(t)\geq\chi_a(t)$ for
	all $t\geq0$.\/}
\end{quote}
We note that $\dot\chi(0)=\cos\theta$, where $\theta$ is the angle between
$\dot\gamma(0)$, and $\dot\beta(0)$.  Let $p\in\B_\gamma(T) \cap
\B_{-\gamma}(T)$, and let $\beta$ be the geodesic from $\gamma(0)$ to $p$.
Let $\theta$ be the angle between $\dot\beta(0)$ and $\dot\gamma(0)$, and
assume first that $\theta\geq\pi/2$.  Then by the comparison principle
above, we have
\[
	\frac1a\, \log\left( e^{ar} \sin^2(\theta/2) + e^{-ar}
	\cos^2(\theta/2) \right) \leq f_\gamma(p) \leq T,
\]
where $r=\dist(p,\gamma(0))$.  Consequently, we obtain
$r\leq T + a^{-1}\log2$.
Now, if $\theta<\pi/2$, then $\pi-\theta$, the angle between $\dot\beta(0)$
and $-\dot\gamma(0)$, is $>\pi/2$, and a similar estimate using
$f_{-\gamma}$ gives $r\leq T + a^{-1}\log2$ again.  The lemma follows.
\end{pf}

\begin{lemma}	\label{constant}
Let $\gamma$ and $\beta$ be geodesics in $\M$, such that
$\beta(-\infty)=\gamma(-\infty)$ and
$\beta(+\infty)\ne\gamma(+\infty)$.  Then for some $d\in\R$:
\begin{align}
	\label{minus}
	\lim_{t\to\infty} \bigl(f_{-\gamma}
	- f_{\gamma}\bigr) \compose\beta(t) &=d, \\
	\label{plus}
	\lim_{t\to-\infty} \bigl(f_{-\gamma}
	+ f_{\gamma}\bigr) \compose\beta(t) &=0.
\end{align}
\end{lemma}
\begin{pf}
Let $\chi = f_{-\gamma}\compose\beta$, and $\chib=f_{\gamma}\compose\beta$.
We have $\dot\chi=1$, while $\dot\chib\leq1$.
Thus, $\chi-\chib$ is non-decreasing
and to prove~\eqref{minus} it remains to show that it is bounded above.
Let $p$ be the unique point where $\beta$ intersects
$S_{-\gamma}(0)$.  Then, for $t$ large enough $\chi(t) = \dist(p,\beta(t))$.
Also, $\chib\to+\infty$ as $t\to\pm\infty$, hence there is $t_0\in\R$ where
$\chib$ has its minimum, and at this point $\dot\beta(t_0)$ is tangent to
$\S_\gamma(t_1)$, where $t_1=f_\gamma\compose\beta(t_0)$.  For each
$t\in\R$, let $\alpha_t$ be the unique geodesic from $\gamma(+\infty)$ to
$\beta(t)$, and let $q_t$ be the unique point where $\alpha_t$ intersects
$S_\gamma(t_1)$.  Then, for $t$ large enough,
$\chib(t)=\dist(q_t,\beta(t))-t_1$, and from the
triangle inequality, we obtain that $\chi(t)-\chib(t)
\leq \dist(p,q_t)+t_1$.
By Theorem~4.9 in~\cite{heintze}, $q_t$ lies in a compact
set, hence $\chi-\chib$ is bounded above, and~\eqref{minus} follows.
Now, note that for any $p\in M$, we have by the triangle inequality
\[
	\bigl[\dist\bigl(p,\gamma(-t)\bigr) - t\bigr] +
	\bigl[\dist\bigl(p,\gamma(t)\bigr) - t\bigr] \geq 0,
\]
hence $\chi+\chib\geq0$.  Also $\dot\chib\geq-1$, thus $\chi+\chib$ is
non-decreasing, and we deduce that $\lim_{t\to-\infty}(\chi+\chib)\geq0$.
It remains to show that $\lim_{t\to-\infty}(\chi+\chib)\leq0$.
We may assume that $\beta$ is parameterized so that $f_{-\gamma} =
f_{-\beta}$.  For any $s<t\leq0$, we have by the triangle inequality:
\[
	\bigl[\dist\bigl(\beta(t),\gamma(-s)\bigr) - s\bigr] +
	\bigl[\dist\bigl(\beta(t),\gamma(s)\bigr) - s\bigr] \leq
	2 \dist(\beta(t),\gamma(t)),
\]
hence, taking the limit $s\to\infty$, we find $\chi(t)+\chib(t) \leq
2\dist(\beta(t),\gamma(t))$. Since by Lemma~\ref{asymptotic},
$\dist(\beta(t),\gamma(t))\to0$ as $t\to-\infty$, Equation~\eqref{plus}
follows.
\end{pf}

\noindent{\em Remark\/}.
{}From the proof it follows that since
$\chi-\chib$ is non-decreasing, the following inequality holds:
\begin{equation}	\label{gamma:i}
	\chi(0) - \chib(0) \leq
	\chi(t) - \chib(t) \leq d,
	\quad \forall t\geq0.
\end{equation}

\vspace{1em}

\begin{lemma}	\label{monotone}
Let $\gamma$ be a geodesic in $\M$.  Then,
for any $t_0\in\R$, and any $T\geq (2a)^{-1}\log2$, we have
\[
	g\bigl(\nabla f_\gamma(p), \nabla f_{-\gamma}(p)\bigr) > 0, \quad
	\forall p\in \S_\gamma(-t_0 + T) \bs \B_{-\gamma}(t_0 + T).
\]
\end{lemma}
\begin{pf}
It is clearly enough to prove the lemma for $T=t_0=(2a)^{-1}\log2$.
Let $p\in\S_\gamma(0)\bs\B_{-\gamma}(a^{-1}\log2)$, and
let $\alpha$ be the unique geodesic from $\gamma(-\infty)$ through $p$
parameterized so that $\alpha(0)\in\S_{-\gamma}(0)$.  Let
$\chib=f_\gamma\compose\alpha$, then $\dot\chib=g(\nabla f_\gamma,\nabla
f_{-\gamma})$.  Now, $\chib$ is convex hence has at most two zeros $t_1\leq
t_2$, and $\dot\chib(t_1)\leq0$.  Let $q=\alpha(t_1)$.
We will show that $t_1\leq a^{-1}\log2$.  That proves the lemma since then
$f_{-\gamma}(q) = t_1 \leq a^{-1}\log2 < f_{-\gamma}(p)$, hence
$p=\alpha(t_2)$ where by the convexity of $\chib$, we have
\[
	g\bigl(\nabla f_\gamma(p),\nabla f_{-\gamma}(p)\bigr) =
	\dot\chib(t_2) > 0.
\]
Let $\beta$ be the  unique geodesic from $\gamma(\infty)$ through $q$
parameterized to that $\beta(0)=q$.  Let $\chi=f_{-\gamma}\compose\beta$,
then $\chi(0)=t_1$.  Thus,
if $\theta$ is the angle between $\beta$ and $\nabla f_{-\gamma}$,
we have $\cos\theta = \dot\chi(0) = \dot\chib(t_1) \leq 0$.  It
follows that $\theta\geq\pi/2$, or equivalently $\sin^2(\theta/2)\geq1/2$.
Since $f_\gamma\compose\beta(0)=0$,
we have $f_\gamma\compose\beta(t)=t$, and
therefore $\chi(t) + t = (f_{-\gamma} + f_\gamma)\compose\beta(t)$.  By
Lemma~\ref{constant}, Equation~\eqref{minus}, we have
\[
	\lim_{t\to-\infty} (\chi(t) + t) = 0.
\]
However, by the comparison principle in Lemma~\ref{bound}, we have
\[
	\chi(t) + t \geq \chi(0) + a^{-1}\, \log\sin^2(\theta/2) +
	a^{-1}\, \log\bigl(1 + e^{at} \cot^2(\theta/2)\bigr).
\]
It follows that $t_1 = \chi(0) \leq -a^{-1} \log\sin^2(\theta/2) \leq
a^{-1}\log2$.
\end{pf}

In order to obtain the next lemma,
which is a key ingredient in the proof of Theorem~\ref{main},
we now assume that $\M$ is a classical
Riemannian globally symmetric space of rank one and of non-compact type.
Thus, $\M$ is $\Ha^\ell_{\K}$, where $\ell\geq2$,
and $\K$ is either $\R$,
$\Co$ or the quaternions $\Ha$.  It is well-known that $\M$ is
simply connected, and when scaled appropriately
has sectional curvatures between $-4$ and $-1$, see~\cite{helgason}.
Thus all of the above considerations apply.
We note that when $\K=\R$, $m=\dim M = \ell$, while when $\K=\Co$,
$m=2\ell$, and when $\K= \Ha$, $m=4\ell$.
Let $\frak S$ denote the sum over cyclic permutations of the indices
$\{1,2,3\}$.

\begin{lemma}	\label{lemma:Qequiv}
Let $\M=\Ha^\ell_{\K}$, where $\ell\geq2$, and
$\K$ is either $\R$, $\Co$ or $\Ha$.  Then
there is an analytic coordinate system $\phi=(u,v)$ on $\M$,
with $u=f_{-\gamma}$, such that the metric $g$ is given in this coordinate
system by the following line elements.\\
When $\K = \R${\em:}
\begin{equation}	\label{realhyperbolic}
	ds^2 = du^2 + e^{2u} \sum_{k=1}^{\ell-1} dv_k^2.
\end{equation}
When $\K = \Co${\em:}
\begin{equation}	\label{complexhyperbolic}
	ds^2 = du^2 + e^{4u} \biggl( dv_1 + \sum_{k=1}^{\ell-1}
	\bigl(v_{2k}\, dv_{2k+1} - v_{2k+1}\, dv_{2k}\bigr)
	\biggl)^2 + \, e^{2u} \sum_{k=2}^{2\ell-1} dv_{k}^2.
\end{equation}
When $\K =  \Ha${\em:}
\begin{equation}	\label{quaternionhyperbolic}
	\begin{split}
	ds^2 & = du^2 \\
	{} &+ e^{4u} \frak S
	\biggl( \! dv_1 \! + \! \sum_{k=1}^{\ell-1} \bigl(
	v_{4k} dv_{4k+1} \! - v_{4k+1} dv_{4k} \! -
	v_{4k+2} dv_{4k+3} \! + v_{4k+3}
	dv_{4k+2} \bigr) \! \biggl)^2\\
	{} &+ e^{2u} \sum_{k=4}^{4\ell-1} dv_k^2.
	\end{split}
\end{equation}
In this coordinate system, the following holds.
\begin{itemize}
\item[(i)]
Let $R>0$, and let $\gamma'$ be a geodesic such that
$\gamma'(-\infty)=\gamma(-\infty)$.
Then there exists $c\geq1$ such that for all $t_0\geq0$, and all
$p\in B_R(\gamma'(t_0))$, there holds:
\begin{equation}	\label{M}
	\frac1c \, Q_{\gamma'(t_0)}(\xi) \leq Q_{p}(\xi)
	\leq c \, Q_{\gamma'(t_0)}(\xi),
	\quad \forall\xi\in\R^{m-1}.
\end{equation}
\item[(ii)]
For all $t,t'\in\R$,
$\S_{-\gamma}(t)\cap\B_{\gamma'}(t')$ is star-shaped in this coordinates
with respect to its `center', the unique point where $\gamma'$ intersects
$\S_{-\gamma}(t)$.
\end{itemize}
\end{lemma}
\begin{pf}
The construction of the coordinate system $\phi$, and the derivation of
Equations~\eqref{realhyperbolic}--\eqref{quaternionhyperbolic} is,
although straightforward, quite
tedious.  We defer it to the appendix.
We turn to the proof of (i).  Note first that
if $p=\gamma'(t)$, then $\norm{t-t_0}\leq\dist(p,\gamma'(t_0))\leq R$,
and hence~\eqref{M} holds with $c=e^{4R}$.
Thus it suffices to prove~\eqref{M} with $\gamma'(t_0)$ replaced by
$\gamma'(t)$, where $t=f_{-\gamma}(p)$, i.e.\
\begin{equation}	\label{M'}
	\frac1{c'} \, Q_{\gamma'(t)}(\xi) \leq Q_{p}(\xi)
	\leq c' \, Q_{\gamma'(t)}(\xi),
	\quad \forall\xi\in\R^{m-1}.
\end{equation}
This will imply~\eqref{M} with $c=e^{4R}c'$.  To prove~\eqref{M'}, consider
first the case $\K=\R$.  In view of Equation~\eqref{realhyperbolic},
$Q_p(\xi)$ is constant on the horospheres, hence there is nothing to prove.
Now consider the case $\K=\Co$.
We may assume that $\phi\compose\gamma'(t) = (t,w)$, where
$w=(w_1,\dots,w_{m-1})\in\R^{m-1}$ is constant.
Note that the map $\tau_{-t}\colon\M\to\M$ given by:
\[
	\phi\compose\tau_{-t}\compose\phi^{-1}(u,v) =
	\bigl( u-t, e^{2t} v_1, e^{t} v_2, \dots, e^{t} v_{m-1} \bigr),
\]
is an isometry, and $\tau_{-t}(\gamma'(t))=\gamma'(0)$.
Also, if $p\in B_R(\gamma'(t_0))\cap\S_{-\gamma}(t)$, then
$\dist(p,\gamma'(t)) \leq \dist(p,\gamma'(t_0)) +
\dist(\gamma'(t_0),\gamma'(t)) \leq 2R$, hence
$p\in B_{2R}(\gamma'(t))\cap\S_{-\gamma}(t)$.
Thus, $\tau_{-t}(p) \in B_{2R}(\gamma'(0))$.  In particular, there is a
constant $C>0$, independent of $t$, such that
$e^{t}\norm{v_k-w_k}\leq C$ for $2\leq k\leq m-1$.
Using~\eqref{complexhyperbolic}, we can therefore estimate:
\begin{align*}
	Q_p(\xi) &=
	e^{4t}\biggl( \xi_1 + \sum_{k=1}^{\ell-1} \bigl(
	v_{2k} \xi_{2k+1} - v_{2k+1} \xi_{2k} \bigr) \biggl)^2
	+ \, e^{2t} \sum_{k=2}^{m-1} \xi_k^2 \\
	&\leq
	2 e^{4t}\biggl( \xi_1 + \sum_{k=1}^{\ell-1} \bigl(
	w_{2k} \xi_{2k+1} - w_{2k+1} \xi_{2k} \bigr) \biggl)^2
	+ \, e^{2t} \sum_{k=2}^{m-1} \xi_k^2 \\
	& \qquad {} + 4 e^{4t} \sum_{k=1}^{\ell-1} \biggl(
	\bigl( v_{2k} - w_{2k} \bigr)^2 \xi_{2k+1}^2 +
	\bigl( v_{2k+1} - w_{2k+1} \bigr)^2 \xi_{2k}^2 \biggl) \\
	& \leq c'\, Q_{\gamma'(t)}(\xi),
\end{align*}
with $c' = \max\{ 2, 4C^2+1 \}$.  The other inequality follows by
interchanging $v$ and $w$.  The case $\K=\Ha$ is proved similarly.
It remains to prove (ii).  Observe that there is a
subgroup $N$ of the group of isometries of $\M$ which leaves each horosphere
$\S_{-\gamma}(t)$ invariant, is transitive on each $\S_{-\gamma}(t)$, and is
linear in $v$.  For instance, when $\K=\Co$, these isometries
$\tau\in N$ are given by
\begin{multline*}
	\phi\compose\tau\compose\phi^{-1}(u,v) = \\
	\bigl( u, \> v_1 + w_1 - \sum_{k=1}^{\ell-1} (w_{2k} v_{2k+1}
	- w_{2k+1} v_{2k}), \> v_2 + w_2, \, \dots,
	\> v_{m-1} + w_{m-1} \bigr),
\end{multline*}
where $w = (w_1,\dots,w_{m-1})\in\R^{m-1}$ is an arbitrary constant.
The case $\K=\Ha$ is similar, and the case $\K=\R$ is trivial.
Thus, we can find such an isometry $\tau\in N$
which maps $\gamma'$ onto $\gamma$, and
leaves $f_{-\gamma}$ invariant.
Since $\tau$ maps `lines' $\phi^{-1}(u,v + tw)$ to `lines'
$\phi^{-1}(u,tv)$, we deduce that it suffices to check (ii) for $\gamma$
only.  The Busemann function $f_{\gamma}$ is computed in the appendix.
Set $\ub=f_{\gamma}\compose\phi^{-1}$, then we have
\begin{equation}	\label{busemann}
	e^{2\ub} =
	\begin{cases}
		\biggl(e^{-u} + e^u \sum_{k=1}^{m-1} v_k^2\biggr)^2 &
		\text{when $\K=\R$}, \\
		\biggl(e^{-u} + e^u \sum_{k=2}^{m-1} v_k^2\biggr)^2 +
		4 e^{2u} v_1^2 & \text{when $\K=\Co$}, \\
		\biggl(e^{-u} + e^u \sum_{k=4}^{m-1} v_k^2\biggr)^2 +
		4 e^{2u} \bigl(v_1^2 + v_2^2 + v_3^2\bigr) &
		\text{when $\K=\Ha$}.
	\end{cases}
\end{equation}
Now it is clear from Equation~\eqref{busemann} that $\ub$ is monotonically
increasing with respect to $\norm{v}$, hence (ii) follows.
\end{pf}

\noindent{\em Remark\/}.  In the proof of Theorem~\ref{main}, we only use
(i) and (ii) of Lemma~\ref{lemma:Qequiv}.  It would be interesting to see
whether these generalize to other simply connected manifolds with
pinched negative curvature.
We note that, according to the proof of Theorem~\ref{main},
it suffices to have (i) and (ii)
for sufficiently large $t$'s.

\vspace{1em}

Let $\Omega\subset\R^n$, $n\geq2$, be a bounded domain
with $\D\Omega$ of class $C^{2,\alpha}$, and let $\Sigma$
be a closed smooth
submanifold of $\Omega$ of co-dimension at least $2$, possibly
with $\D \Sigma\ne\emptyset$.

\begin{defn}
Let $\gamma$ be a geodesic in $\M$.
We say that a harmonic map $\vp\in C^\infty(\Omega\bs\Sigma;M)\cap
C^{2,\alpha}(\overline\Omega\bs\Sigma;M)$ is a {\em
$\Sigma$-singular map into $\gamma$}
if
\begin{itemize}
	\item[(i)] $\vp(\Omega\bs\Sigma)\subset\gamma(\R)$
	\item[(ii)] $\vp(x)\to\gamma(+\infty)$ as $x\to\Sigma$
	\item[(iii)] There is a constant $\delta>0$
	such that $\norm{d\vp(x)}^2 \geq \delta \dist(x,\Sigma)^{-2}$
	in a neighborhood of $\Sigma$.
\end{itemize}
\end{defn}

Let $\mu$ be any positive measure
on $\Sigma$, and let $\Gamma$ be the fundamental solution in $\Omega$,
then the convolution $u=\mu*\Gamma$ is a
harmonic function on $\Omega\bs\Sigma$
which tends to infinity on $\Sigma$,
and hence $\vp=\gamma\compose u$ satisfies (i) and (ii).
Conversely, if $\vp$ is a $\Sigma$-singular map into $\gamma$, then there is
a harmonic function $u$ on $\Omega\bs\Sigma$ such that
$\vp=\gamma\compose u$,
and $u=\mu*\Gamma-u'$ for some positive measure $\mu$ on $\Sigma$, and
some smooth harmonic function $u'$ on $\Omega$.
Since $\vp=\gamma\compose u$ implies $\norm{d\vp}^2 = \norm{\nabla u}^2$,
condition (iii) can be obtained for example if the measure
$\mu$ is bounded below by a positive constant $\delta$ times
the surface measure of $\Sigma$.
Note that if $u$ and $u'$ are harmonic functions on $\Omega\bs\Sigma$, and
$u-u'$ is a smooth harmonic function on $\Omega$, then
$\gamma\compose u$ and $\gamma\compose u'$ are asymptotic.
Thus, if $\vp=\gamma\compose u$ is a $\Sigma$-singular map into $\gamma$,
we can always assume without loss of generality
that $\vp$ maps $\D\Omega$ to $\gamma(0)$, for otherwise we
can add to $u$ a smooth harmonic function $u'$ so that $u+u'=0$ on
$\D\Omega$.  In particular, we may assume that $u>0$ in
$\Omega\bs\Sigma$.

Let $L^\infty(\Omega\bs\Sigma)$ be the space of measurable
functions on $\Omega\bs\Sigma$ which are essentially bounded.
In analogy with geodesics we define:

\begin{defn}
Let $\vp,\vp'\colon\Omega\bs\Sigma\to\M$ be harmonic maps, and let
$\Sigma'\subset\Sigma$.
We say that $\vp$ and $\vp'$ are {\em asymptotic near\/} $\Sigma'$
if there is a neighborhood $\Omega'$ of $\Sigma'$ such that
$\dist(\vp,\vp')\in L^\infty(\Omega'\bs\Sigma')$.  If\/
$\dist(\vp,\vp') \in L^\infty(\Omega\bs\Sigma)$,
we say they are {\em asymptotic\/}.
\end{defn}

We will also use the following two elementary lemmas.
The first is an integral estimate for singular harmonic functions.
The second is a simple maximum principle.
For the sake of completeness, we give the proofs,
although they are quite standard.

\begin{lemma}	\label{lemma:integral:estimate}
Let $u\in C^\infty(\Omega\bs\Sigma)$ satisfy $\Delta u=0$ in
$\Omega\bs\Sigma$, $u>0$ in $\Omega\bs\Sigma$, and
$u(x)\to\infty$ as $x\to\Sigma$.
Let $\D^s=\{x\in\Omega:\> \dist(x,\Sigma)=s\}$.  Then, we have:
\begin{equation}	\label{integral:limit}
	\lim_{s\to0} \int_{\D^s} u =0.
\end{equation}
\end{lemma}
\begin{pf}
First observe that for $s>0$ small enough, $\D^s$ is a smooth
compact surface.  Let $\chi(s)= \int_{\D^s} u\geq0$, then we compute:
\begin{equation}	\label{mean:curvature}
	\frac{d\chi}{ds} = \int_{\D^s} \, \D_n u + \int_{\D^s}
	u \, h,
\end{equation}
where $\D_n u$ is the derivative of $u$ along the outward unit normal to
$\D^s$, and $h$ is the mean curvature of $\D^s$.  Note that, since $u$ is
harmonic, the first integral is independent of $s>0$, and define the {\em
charge\/} of $u$:
\[
	e_0 = -\int_{\D^s} \, \D_n u.
\]
Furthermore, since $\Sigma$ is of co-dimension $k\geq2$, we have
on $\D^s$ for $s>0$ small enough:
\[
	h \geq \frac{k-1}{s} - c_1 \geq \frac1s - c_1,
\]
for some $c_1\geq0$.  To see this, note that if $X$ is the field of unit
normals to the surfaces $\D^s$, then as $s\to0$ the dominant part in
$h=\Div X$ comes from the divergence of $X$ in the $k$-planes
normal to $\Sigma$ which is $(k-1)/s$.
Thus, from~\eqref{mean:curvature} we obtain the differential inequality:
\begin{equation}	\label{chi}
	\frac{d\chi}{ds} \geq -e_0 + \left(\frac1s - c_1\right) \chi.
\end{equation}
Equation~\eqref{integral:limit} follows from~\eqref{chi}.
To see this, note first that Inequality~\eqref{chi} implies that $\chi(s)$
is increasing for $s>0$ small enough.  Thus, we have $\chi(s) \leq
\chi_0=\chi(s_0)$ for $0<s\leq s_0$.  Now, introduce $t=\log s$.
Then Inequality~\eqref{chi} becomes:
\[
	\frac{d}{dt}\biggl(e^{-t}\chi + (e_0+c_1\chi_0)\, t\biggr) \geq 0,
\]
and that yields $\chi\to0$ as $t\to-\infty$.
\end{pf}

\begin{lemma}	\label{uniqueness}
Let $u\in C^2(\Omega\bs\Sigma) \cap C^1(\overline\Omega\bs\Sigma)$
satisfy $\Delta u\geq0$ in
$\Omega\bs\Sigma$, $0\leq u\leq 1$ in $\Omega\bs\Sigma$,
and $u|\D\Omega=0$.  Then $u=0$.
\end{lemma}
\begin{pf}
For any function
$\chi\in C^{0,1}_0(\R^n\bs\Sigma)$, with $0\leq \chi\leq 1$, we have:
\[
	0 \geq -\int_{\Omega} \chi^2 u \, \Delta u = \int_{\Omega} \chi^2
	\norm{\nabla u}^2 + 2 \int_{\Omega} \chi u \nabla\chi\cdot \nabla u.
\]
It follows that
\[
	\int_{\Omega} \chi^2 \norm{\nabla u}^2 \leq 2 \left(\int_{\Omega}
	\chi^2 \norm{\nabla u}^2 \right)^{1/2} \left(\int_{\Omega}
	\norm{\nabla \chi}^2\right)^{1/2}.
\]
and therefore:
\begin{equation}	\label{estimate}
	\int_{\Omega} \chi^2 \norm{\nabla u}^2 \leq 4 \int_{\Omega}
	\norm{\nabla \chi}^2.
\end{equation}
Now, let $r(x) = \dist(x,\Sigma)$, and for $\epsilon>0$ small enough, define
\[
	\chi_\epsilon =
	\begin{cases}
		2 - \log r/\log \epsilon & \text{if $\epsilon^2 \leq r \leq
		\epsilon$} \\
		0 & \text{if $r\leq \epsilon^2$}\\
		1 & \text{if $r\geq\epsilon$}.
	\end{cases}
\]
Then $\chi_\epsilon \in C^{0,1}_0(\R^n\bs\Sigma)$,
$0\leq\chi_\epsilon\leq1$, and
\[
	\nabla\chi_\epsilon =
	\begin{cases}
		- \disp\frac{\nabla r}{r\log\epsilon} &
		\text{if $\epsilon^2\leq
		r\leq \epsilon$}\\
		0 & \text{otherwise}.
	\end{cases}
\]
Let $a(s)$ be the $(n-1)$-Hausdorff measure of $\D^s=\{x\in\R^n:
\dist(x,\Sigma)=s\}$.  Since $\Sigma$ is of co-dimension $\geq2$,
an argument similar to the one used in Lemma~\ref{lemma:integral:estimate}
shows that $a(s) \leq C s$.  Thus, by the co-area formula, see~\cite{simon},
we have:
\[
	\int_{\R^n} \norm{\nabla \chi_\epsilon}^2 =
	\frac{1}{\bigl(\log\epsilon\bigr)^2} \int_{\epsilon^2}^\epsilon
	a(s)\, \frac{ds}{s^2} \leq \frac{C}{\norm{\log\epsilon}}.
\]
Substituting this into Inequality~\eqref{estimate}, we obtain:
\[
	\int_{\Omega} \chi_\epsilon^2 \norm{\nabla u}^2 \leq
	\frac{C}{\norm{\log\epsilon}}.
\]
Now, if we let $\epsilon\to0$, the left hand side tends to $\int_{\Omega}
\norm{\nabla u}^2$, while the right
hand side tends to zero.  Consequently, $u$
is constant, and since $u=0$ on $\D\Omega$, we conclude that $u=0$.
\end{pf}

\section{The Case $N=1$} \label{single}

In this section we prove the following proposition, which is the case $N=1$
in Theorem~\ref{main}.

\begin{prop}	\label{prop:single}
Let $\gamma$ be a geodesic in $\M$, let $\vp_0$ be a $\Sigma$-singular map
into $\gamma$, and let $\psi\in C^{2,\alpha}(\D\Omega;M)$.  Then
there exists a unique
harmonic map $\vp\in C^\infty(\Omega\bs\Sigma;M)\cap
C^{2,\alpha}(\overline\Omega\bs\Sigma;M)$
such that $\vp=\psi$ on
$\D\Omega$, and $\vp$ is asymptotic to $\vp_0$.
\end{prop}
\noindent{\em Proof\/}.
We begin with the uniqueness.  Suppose that $\vp$ and $\vp'$ are harmonic
maps which are asymptotic to $\vp_0$ and agree with $\psi$ on
$\D\Omega$.  Then, $\vp$ and $\vp'$ are asymptotic.  Let $u=\bigl(
\dist(\vp,\vp')\bigr)^2$, then, we have
$u\in C^2(\overline\Omega\bs\Sigma)$, $\Delta u\geq 0$
on $\Omega\bs\Sigma$,
see~\cite{schoen79}, $u$
is bounded, and vanishes on $\D\Omega$.  Thus, in view of
Lemma~\ref{uniqueness}, it follows that $u=0$, hence $\vp=\vp'$.

To prove the existence, we set up a variational principle.  Let
$u=f_{-\gamma}$,
and let $\phi=(u,v)$ be the corresponding coordinate system on $M$
given in Lemma~\ref{lemma:Qequiv}.
Where no confusion arises, we will
identify $\vp$ and its parameterization $\phi\compose\vp=(u,v)$.
Let $\vp_0 = \gamma\compose u_0$,
then $\vp_0 = (u_0, 0)$, $\Delta u_0=0$ on $\Omega\bs\Sigma$, and
we can assume without loss of generality that $u_0=0$ on $\D\Omega$, and
consequently
$u_0>0$ on $\Omega\bs\Sigma$.    In addition, we know that there exists a
constant $\delta>0$ such that
\begin{equation}	\label{norm:u0}
	 \norm{\nabla u_0}^2 = \norm{d\vp_0}^2 \geq \delta r^{-2},
\end{equation}
in a neighborhood of $\Sigma$, where $r(x)=\dist(x,\Sigma)$.
Also, for any $\vp\colon\Omega\bs\Sigma\to M$, $Q_\vp$ is a function from
$\Omega\bs\Sigma$ with values in the positive quadratic forms on $\R^{m-1}$.
Finally, we note that for $\Omega'\subset\Omega$:
\[
	E_{\Omega'}(\vp) = \int_{\Omega'} \left\{
	\norm{\nabla u}^2 + Q_\vp(\nabla v) \right\},
\]
where $Q_\vp(\nabla v) = \sum_{k=1}^n Q_\vp(\nabla_k v)$.

Let $H_{1}(\Omega)$ be the Sobolev space of
functions $u$ such that $u$ and $\nabla u\in L^2(\Omega)$, and
let $H_{1,0}(\Omega)$ be the closure of $C^\infty_0(\Omega)$ in
that space with respect to the norm:
\[
	\Vert u \Vert =
	\left( \int_{\Omega} \left\{ u^2 + \norm{\nabla u}^2 \right\}
	\right)^{1/2}.
\]
We define the weighted Sobolev space $H_{1}^{\vp_0}(\Omega;\R^{m-1})$
to be the space of functions $v\in L^2(\Omega\bs\Sigma;\R^{m-1})$ such that:
\begin{equation}	\label{Qnorm}
	\Vert v \Vert_{\vp_0} =
	\left( \int_{\Omega} \left\{ \norm{v}^2
	+ Q_{\vp_0}(\nabla v) \right\}
	\right)^{1/2}
	< \infty,
\end{equation}
and we define $H_{1,0}^{\vp_0}(\Omega;\R^{m-1})$ to be the closure of
$C^\infty_0(\Omega\bs\Sigma;\R^{m-1})$ in that space with respect to
the norm~\eqref{Qnorm}.  Note that for each $\xi\in\R^{m-1}$,
$Q_{\vp_0(x)}(\xi)\geq Q_{\gamma(0)}(\xi)\geq c_1 \norm{\xi}^2$ for some
$c_1>0$, where the last norm is the Euclidean norm on $\R^{m-1}$.  Thus, we
have the continuous embedding $H_1^{\vp_0}(\Omega;\R^{m-1}) \hookrightarrow
H_1(\Omega;\R^{m-1})$.  In particular, due to the Poincar\'e inequality in
$H_{1,0}(\Omega;\R^{m-1})$:
\[
	\int_\Omega \norm{v}^2 \leq C \int_\Omega \norm{\nabla v}^2,
	\qquad \forall v\in H_{1,0}(\Omega;\R^{m-1}),
\]
we have that the semi-norm:
\[
	\left( \int_\Omega Q_{\vp_0}(\nabla v) \right)^{1/2}
\]
is equivalent on
$H_{1,0}^{\vp_0}(\Omega;\R^{m-1})$ to the full norm~\eqref{Qnorm}.
We also note the following weighted Poincar\'e inequality in
$H_{1,0}^{\vp_0}(\Omega;\R^{m-1})$, compare
with~\cite[Lemma~1]{weinstein92}:

\begin{lemma}	\label{poincare}
For every $v\in H_{1,0}^{\vp_0}(\Omega;\R^{m-1})$, there holds:
\begin{equation}	\label{eq:poincare}
	\int_{\Omega} Q_{\vp_0}(v)\, \norm{\nabla u_0}^2 \leq
	a^{-2} \int_{\Omega} Q_{\vp_0}(\nabla v).
\end{equation}
\end{lemma}
\begin{pf*}{Proof of Lemma~\ref{poincare}}
By a standard density argument, it suffices to
prove~\eqref{eq:poincare} for every $v\in
C^\infty_0(\Omega\bs\Sigma;\R^{m-1})$, and by extending $v$ to be zero
outside $\Omega$, we may assume that $v\in
C^{\infty}_0(\R^n\bs\Sigma;\R^{m-1})$.  Let $S_t = \{x\in\R^{n}:
u_0(x) = t\}$, then for $t$ large enough, by~\eqref{norm:u0},
$S_t$ is a smooth hypersurface, with
interior unit normal $\norm{\nabla u_0}^{-1} \nabla u_0$.
Since for $t$ large enough, we have $v=0$ on $S_t$, we obtain:
\[
	\int_{B_t} \Div\bigl( Q_{\vp_0}(v) \, \nabla u_0\bigr) =
	\int_{S_t} Q_{\vp_0}(v)\, \norm{\nabla u_0} = 0,
\]
where $B_t = \{x\in\R^n: u_0(x)\leq t\}$.  Furthermore,
in view of~\eqref{expansion} Lemma~\ref{lemma:expansion}, we can estimate:
\begin{align*}
	\Div\bigl(Q_{\vp_0}(v)\, \nabla u_0\bigr) &=
	\frac{\D}{\D u} \bigl(Q_{\vp_0}(v)\bigr) \, \norm{\nabla u_0}^2 +
	2 \, Q_{\vp_0}(v,\nabla v) \, \nabla u_0\\
	& \geq  a\, Q_{\vp_0}(v)\, \norm{\nabla u_0}^2 -
	a^{-1}\, Q_{\vp_0}(\nabla v).
\end{align*}
Consequently, we conclude that for $t$ large enough, there holds:
\[
	\int_{B_t} Q_{\vp_0}(v) \, \norm{\nabla u_0}^2 \leq
	a^{-2}\, \int_{B_t} Q_{\vp_0} (\nabla v).
\]
The lemma follows by taking $t\to\infty$.
\end{pf*}

Now, extend $\psi$ to a map $\tilde\psi\in C^\infty(\Omega;M)\cap
C^{2,\alpha}(\overline\Omega;M)$ which
maps a neighborhood of $\Sigma$ to the
point $\gamma(0)$, and write
$\tilde\psi=(\tilde u,\tilde v)$.  Then $(\tilde
u,\tilde v) = (0,0)$ in a neighborhood of $\Sigma$.
Define $\H$ to be the space of maps $\vp=(u,v)\colon\Omega\bs\Sigma\to\M$
such that $\dist(\vp,\vp_0)\in L^\infty(\Omega\bs\Sigma)$,
$u-u_0-\tilde u\in H_{1,0}(\Omega)$, and
$v-\tilde v\in H_{1,0}^{\vp_0}(\Omega;\R^{m-1})$.

For maps $\vp\in\H$, and $\Omega'\subset\Omega$,
we define:
\[
	F_{\Omega'}(\vp) =
	\int_{\Omega'} \left\{ \norm{\nabla(u - u_0)}^2 + Q_\vp(\nabla
	v) \right\},
\]
and we set $F = F_\Omega$.  Note that $F\geq0$, and Lemma~\ref{lemma:equiv}
below implies that $F<\infty$ on $\H$.  We first
show that if $\vp\in\H$ is a
minimizer of $F$, then $\vp$ is a harmonic map on $\Omega\bs\Sigma$, hence
$\vp\in C^\infty(\Omega\bs\Sigma;M)\cap
C^{2,\alpha}(\overline\Omega\bs\Sigma;M)$, and $\vp$ is asymptotic to
$\vp_0$.  Indeed, let $\Omega'\subset\subset\Omega\bs\Sigma$,
then we claim that for any map
$\vp'\in C^1(\Omega';M)$ of finite energy such that $\vp'=\vp$
on $\D\Omega$,
we have $E_{\Omega'}(\vp)\leq E_{\Omega'}(\vp')$.
To see this, note that
\begin{align*}
	E_{\Omega'}(\vp) - F_{\Omega'}(\vp)
	& = \int_{\Omega'} \nabla(2u-u_0) \cdot \nabla u_0 \\
	& = \int_{\Omega'} \Div\bigl( (2u-u_0) \nabla u_0\bigr) \\
	& = \int_{\D\Omega'} (2u-u_0) \, \frac{\D u_0}{\D n} \\[1ex]
	& = E_{\Omega'}(\vp') - F_{\Omega'}(\vp').
\end{align*}
Thus, if instead we had $E_{\Omega'}(\vp')<E_{\Omega'}(\vp)$, then also
$F_{\Omega'}(\vp')<F_{\Omega'}(\vp)$ would hold.
Then the map $\vp''$ defined by:
\[
	\vp''(x) = \begin{cases}
		\vp(x) & \text{if $x\not\in\Omega'$}\\
		\vp'(x) & \text{if $x\in\Omega'$}
	\end{cases}
\]
would satisfy $\vp''\in\H$ and $F(\vp'')<F(\vp)$, in contradiction
to $\vp$ being a minimizer.  Thus
$\vp$ is a critical point of $E_{\Omega'}$, and it follows that
$\vp$ is harmonic.  The interior regularity statement is standard.
Clearly, $\vp=\psi$ on $\D\Omega$, hence the
boundary regularity statement follows.  Therefore, to
prove Proposition~\ref{prop:single}, it suffices to show that
$F$ has a minimizer in $\H$.

For any $R>0$, define the space $\Hc$ of
maps $\vp\in\H$ for which $\dist(\vp,\vp_0)\leq R$ for a.e.\
$x\in\Omega\bs\Sigma$.
We first we show that $F$ has a minimizer on $\Hc$.
For this purpose, we will need the following lemma:

\begin{lemma}	\label{lemma:equiv}
Let $R>0$, then there is $c\geq1$ such that for all $\vp\in\Hc$,
there holds:
\begin{equation}	\label{equivalent}
		\frac1c \, Q_{\vp_0(x)}(\xi)
		\leq Q_{\vp(x)}(\xi) \leq c \,
		Q_{\vp_0(x)}(\xi),\quad
		\forall\xi\in\R^{m-1},\>
		\text{a.e.\ } x\in\Omega\bs\Sigma
\end{equation}
\end{lemma}
\begin{pf*}{Proof of Lemma~\ref{lemma:equiv}}
For every $x\in\Omega\bs\Sigma$ such that
$\dist((\vp(x),\vp_0(x))\leq R$, we have $\vp(x)\in
B_R(\vp_0(x))$, and there is $t\geq0$ such that
$\vp_0(x)=\gamma(t)$.  Thus the lemma follows from (i) in
Lemma~\ref{lemma:Qequiv}.
\end{pf*}

Now, let $\vp_j=(u_j,v_j)\in\Hc$ be a minimizing sequence.
Then $u_j-u_0$ is bounded in
$H_1(\Omega)$, and by passing to a subsequence, we may assume that
$u_j-u_0$
converges weakly and pointwise a.e.\ in $\Omega$ to $u-u_0\in H_1(\Omega)$.
Clearly, we have
\begin{equation}	\label{uest}
	\int_{\Omega} \norm{\nabla(u-u_0)}^2 \leq \liminf \int_{\Omega}
	\norm{\nabla(u_j-u_0)}^2,
\end{equation}
and $u-u_0-\tilde u\in H_{1,0}(\Omega)$.
Now, it follows from~\eqref{equivalent} that $v_j-\tilde v$ is bounded in
$H_{1,0}^{\vp_0}(\Omega;\R^{m-1})$:
\begin{align*}
	\int_{\Omega} Q_{\vp_0}(\nabla v_j - \nabla \tilde v) &\leq
	2 c \int_\Omega Q_{\vp_j}(\nabla v_j) + 2 \int_\Omega
	Q_{\vp_0}(\nabla \tilde v) \\
	& \leq 2c\, (\, \inf_{\Hc} F + c_1) + 2
	\int_\Omega Q_{\vp_0} (\nabla \tilde v),
\end{align*}
for some $c_1\geq0$.  Hence by passing again to a subsequence, we may
assume that $v_j-\tilde v$
converges weakly and pointwise a.e.\ in $\Omega$ to
$v-\tilde v\in H_{1,0}^{\vp_0}(\Omega;\R^{m-1})$.
It follows at once that $\vp=(u,v)\in\H$.
Furthermore, since $\vp_j(x)\to\vp(x)$ for a.e.\ $x\in\Omega$, we obtain
that $\vp\in\Hc$.  Thus, by Lemma~\ref{lemma:equiv},
the quadratic form $Q_\vp$ is uniformly
equivalent to $Q_{\vp_0}$, hence we have:
\begin{equation}	\label{weak}
	\int_{\Omega} Q_\vp (\nabla v) =
	\lim \int_{\Omega} Q_\vp (\nabla v, \nabla v_j),
\end{equation}
where on the right hand side, we have used $Q_\vp$ also for the symmetric
bilinear form associated with $Q_\vp$.
Define
\[
	\chi_j = \begin{cases}
		\disp\frac{Q_\vp(\nabla v_j)}{Q_{\vp_j}(\nabla
		v_j)} &\text{if $\nabla v_j\ne 0$} \\[.5ex]
		1 & \text{otherwise},
		\end{cases}
\]
then, by~\eqref{equivalent}, we have $\chi_j\leq c$, and furthermore,
we claim, $\chi_j\to1$ pointwise a.e.\ in $\Omega$.
To see this, suppose that $p_j=\vp_j(x)\to p=\vp(x)$,
then we will show that $\chi_j(x)\to1$.
Since $Q_p$ is continuous on $M$, we have
$Q_{p_j}(\xi) \to Q_{p}(\xi)$
uniformly for $\xi\in \sphere^{m-2}$, the unit sphere in $\R^{m-1}$.
Let $j'$ be the subsequence of $j$'s for which
$\nabla v_{j'}(x) \ne 0$.
Clearly, it is sufficient to prove that $\chi_{j'}(x)\to1$.
For every $j'$ and
for each $1\leq k\leq n$, there are $\lambda_{j'}^k\geq0$
and $\xi_{j'}^k\in \sphere^{m-2}$, such that $\nabla_k v_{j'}(x) =
\lambda_{j'}^k\, \xi_{j'}^k$, and $\lambda_{j'} = \max_k\lambda_{j'}^k>0$.
We see that at $x$, we have:
\[
	Q_{p_{j'}}(\nabla v_{j'}) =
	\sum_{k=1}^n \left(\lambda_{j'}^k\right)^2
	Q_{p_{j'}}(\xi_{j'}^k)
	\geq \left(\lambda_{j'}\right)^2 \inf_{\sphere^{m-2}} Q_{p_{j'}}
	\geq \frac12\,
	\left(\lambda_{j'}\right)^2 \inf_{\sphere^{m-2}} Q_p,
\]
for $j'$ large enough.
Consequently, if we set
$c_2=2 \, \bigl(\inf_{\sphere^{m-2}} Q_{p}\bigr)^{-1}$,
we can conclude that:
\[
	\norm{\chi_{j'}(x) - 1} \leq c_2\sum_{k=1}^n
	\left(\frac{\lambda_{j'}^k}{\lambda_{j'}}\right)^2
	\norm{ Q_{p}(\xi_{j'}^k) - Q_{p_{j'}}(\xi_{j'}^k) }
	\leq n c_2 \sup_{\sphere^{m-2}} \norm{ Q_p - Q_{p_{j'}} },
\]
which tends to zero when $j'\to\infty$.  We can now estimate the right hand
side of~\eqref{weak}:
\begin{equation}	\label{vest}
	\begin{split}
	\int_{\Omega} Q_{\vp}(\nabla v,\nabla v_j) &\leq
	\int_{\Omega} \biggl( Q_{\vp}(\nabla v) \biggr)^{1/2}
	\biggl( Q_{\vp}(\nabla v_j) \biggr)^{1/2} \\
	&\leq \left( \int_{\Omega} \chi_j \, Q_{\vp}(\nabla v)
	\right)^{1/2} \left( \int_{\Omega} Q_{\vp_j}(\nabla v_j)
	\right)^{1/2}.
	\end{split}
\end{equation}
The first factor on the right hand side of~\eqref{vest} has a
limit by the Dominated Convergence Theorem:
\begin{equation}	\label{limit}
	\lim \int_{\Omega} \chi_j \, Q_\vp(\nabla v)
	= \int_{\Omega} Q_\vp(\nabla v).
\end{equation}
Combining~\eqref{limit} with~\eqref{weak} and~\eqref{vest}, we obtain:
\[
	\int_{\Omega} Q_\vp (\nabla v) \leq
	\left( \int_{\Omega} Q_\vp (\nabla v) \right)^{1/2}
	\left( \liminf \int_{\Omega}
	Q_{\vp_j}(\nabla v_j) \right)^{1/2},
\]
from which it follows that:
\begin{equation}	\label{vest1}
	\int_{\Omega} Q_\vp (\nabla v)
	\leq \liminf \int_{\Omega}
	Q_{\vp_j}(\nabla v_j).
\end{equation}
In view of~\eqref{uest} and~\eqref{vest1}, we conclude that:
\[
	F(\vp) \leq \liminf \int_{\Omega} \norm{\nabla{(u_j-u_0)}}^2 +
	\liminf \int_{\Omega}
	Q_{\vp_j}(\nabla v_j) \leq \lim F(\vp_j)
	= \inf_{\Hc} F,
\]
Hence $F(\vp) = \inf_{\Hc} F$, and $\vp$ is a minimizer of $F$ on $\Hc$.

It remains to prove that for some $R>0$ large enough
$\inf_{\H} F = \inf_{\Hc} F$.
Clearly, it suffices to prove that $\inf_{\Hc} F \leq \inf_{\H} F$.
This follows from the following lemma.

\begin{lemma}	\label{lemma:inf}
There is a constant $R>0$ such that for every $\epsilon>0$, and every
$\vp\in\H$, there is $\vp'\in\Hc$ such that
$F(\vp')\leq F(\vp) + \epsilon$.
\end{lemma}
\begin{pf*}{Proof of Lemma~\ref{lemma:inf}}
Let $\H^*$ be the space of maps $\vp=(u,v)\in
\H$ such that $v=0$ in a neighborhood of $\Sigma$.
The proof of the lemma will be divided into two steps:
we will prove that (i) {\em given $\epsilon>0$ and $\vp\in\H$,
there is $\vp'\in\H^*$ such
that $F(\vp')\leq F(\vp) + \epsilon$\/}; and (ii)
{\em a constant $R>0$ exists such that given
any $\vp\in\H^*$, there is $\vp'\in\Hc$ such that $F(\vp') \leq
F(\vp)$\/}.  The lemma immediately follows from (i) and (ii).
To prove (i), we use the function
$\chi_\epsilon$ introduced in the proof of
Lemma~\ref{uniqueness}.  Note that, in view of~\eqref{norm:u0}, we have
\[
	\norm{\nabla\chi_\epsilon}^2 =
	\dist(x,\Sigma)^{-2}\norm{\log\epsilon}^{-2} \leq
	\delta^{-1} \norm{\log\epsilon}^{-2} \norm{\nabla u_0}^2.
\]
Let $\vp=(u,v)\in\H$, define $v_\epsilon = \chi_\epsilon v$, and
$\vp_\epsilon = (u,v_\epsilon)$, then
$v_\epsilon - \tilde v \in H_{1,0}^{\vp_0}(\Omega;\R^{m-1})$, and
$\vp_\epsilon\in\H^*$.  Furthermore, we have
$\vp\in\H_{R'}$ for some $R'>0$.  Let $x\in\Omega\bs\Sigma$ be such that
$\dist(\vp(x),\vp_0(x))\leq R'$, then $u(x) \leq u_0(x) + R'$.
Consequently, we have $\vp(x) \in \B_{-\gamma}(u_0(x)+R')$.
Similarly, $\vp(x)\in\B_\gamma(-u_0(x)+R')$.
By (ii) in Lemma~\ref{lemma:Qequiv}, this implies that
$\vp_\epsilon(x)\in\B_\gamma(-u_0(x)+R')$, and therefore
$\vp_\epsilon(x)\in\B_{-\gamma}(u_0(x)+R')\cap\B_{\gamma}(-u_0(x)+R')$.
It follows by Lemma~\ref{bound} that $\vp\in\H_{R''}$ where
$R''=R'+a^{-1}\log 2$.  Now, we have
\begin{equation}	\label{Qe}
	\int_\Omega Q_{\vp_\epsilon}(\nabla v_\epsilon) =
	\int_\Omega\chi_\epsilon^2\,
	Q_{\vp_\epsilon}(\nabla v) + I_\epsilon,
\end{equation}
where
\[
	\norm{I_\epsilon}
	\leq \int_\Omega  Q_{\vp_\epsilon} (v^2)\, \norm{\nabla
	\chi_\epsilon}^2 + 2 \left( \int_\Omega Q_{\vp_\epsilon}(v^2)\,
	\norm{\nabla\chi_\epsilon}^2 \right)^{1/2}
	\left( \int_\Omega \chi_\epsilon^2\,
	Q_{\vp_\epsilon}(\nabla v) \right)^{1/2}.
\]
By Lemma~\ref{lemma:equiv}, we have
\[
	\int_\Omega  Q_{\vp_\epsilon}(v^2)\,\norm{\nabla \chi_\epsilon}^2
	\leq c\, \int_\Omega
	Q_{\vp_0}(v^2)\, \norm{\nabla\chi_\epsilon}^2
	\leq c\, \delta^{-1}\,
	\norm{\log\epsilon}^{-2} \int_{\Omega} Q_{\vp_0}(v^2)
	\, \norm{\nabla u_0}^2,
\]
which, in view of Lemma~\ref{poincare}, tends to zero
as $\epsilon\to0$.  It follows that
$\norm{I_\epsilon}\to0$ as $\epsilon\to0$.  Hence, since
$\chi_\epsilon^2 \, Q_{\vp_\epsilon}(\nabla v)\to Q_{\vp} (\nabla  v)$,
and since
$\chi_\epsilon^2\, Q_{\vp_\epsilon}(\nabla v) \leq c \,
Q_{\vp_0}(\nabla v)$, we deduce from~\eqref{Qe}, by the
Dominated Convergence Theorem, that:
\[
	\int_\Omega Q_{\vp_\epsilon}(\nabla v_\epsilon) \to \int_\Omega
	Q_{\vp}(\nabla v).
\]
Therefore, (i) is obtained.  We now turn to the proof of (ii).  Introduce a
new coordinate system $\phib=(\ub,\vb)\colon M\to\R^m$, where
$\ub=f_\gamma$.  As before, the metric on $\M$ can be written as:
\[
	ds^2 = d\ub\,{}^2 + \Qb_{p} (d\vb).
\]
Let $\vp\in\H^*$, and write $(\ub,\vb)=\phib\compose\vp$.
In $\Omega\bs\Sigma$, we have:
\begin{align*}
	\norm{\nabla(u-u_0)}^2 + Q_\vp(\nabla v) &=
	\norm{\nabla u}^2 + Q_\vp(\nabla v) +  \norm{\nabla u_0}^2 - 2\,
	\nabla u \cdot\nabla u_0 \\
	&=
	\norm{\nabla\ub}^2 + \Qb_\vp(\nabla\vb) + \norm{\nabla u_0}^2 - 2\,
	\nabla u \cdot\nabla u_0 \\
	&=
	\norm{\nabla(\ub+u_0)}^2 + \Qb_\vp(\nabla\vb) - 2\, \Div\bigl(
	(u+\ub) \nabla u_0 \bigr).
\end{align*}
We now integrate this identity over $\Omega^\epsilon = \{ x\in\Omega:
\dist(x,\Sigma)>\epsilon\}$.  If $\epsilon>0$ is small enough,
we can decompose $\D\Omega^\epsilon$ into a disjoint union
$\D\Omega\cup\D^\epsilon$, where $\D^\epsilon = \{x\in\Omega:
\dist(x,\Sigma) = \epsilon\}$.
We observe that if $\epsilon>0$ is small enough, then
$v=0$, and hence $u+\ub=0$ on $\D^\epsilon$.
Consequently, after taking $\epsilon\to0$, we obtain:
\begin{equation}	\label{integral:identity}
	F(\vp) = \int_\Omega \left\{\norm{\nabla(\ub+u_0)}^2 +
	\Qb_\vp(\nabla\vb) \right\} - 2\int_{\D\Omega} (u+\ub) \frac{\D
	u_0}{\D n}.
\end{equation}
Note that the second term can be written as:
\[
	2\int_{\D\Omega} (f_{-\gamma}+f_{\gamma})\compose\psi\,
	\left( \frac{\D u_0}{\D n} \right) ,
\]
which clearly depends only on $\psi$ and $\vp_0$, and hence is
constant in $\H^*$.  Let
\[
	\ub' = -u_0 + \min \{ \, \ub + u_0,\, \Tb\},
\]
where
\[
	\Tb = \sup_{\D\Omega} \ub + 1 = \sup_{\D\Omega} \bigl(
	f_\gamma\compose\psi \bigr)
	+ 1,
\]
and define $\vp'$ by $\phib\compose\vp' = (\ub',\vb)$.
Then $\vp'=\vp$ on $\D\Omega$, hence, in view of~\eqref{integral:identity},
$\vp'\in\H$.  Also,
from~\eqref{integral:identity} and the fact that $\Qb_{\vp}(\nabla \vb)$
is non-decreasing in $\ub$, we have
$F(\vp') \leq F(\vp)$.  Furthermore, the estimate:
\[
	\ub' + u_0 \leq \Tb,
\]
holds throughout $\Omega$.  This estimate can be rewritten
as $f_\gamma\compose\vp' \leq -u_0 + \Tb$, or equivalently as:
\begin{equation}	\label{ub:estimate}
	\vp'(x) \in \B_\gamma(-u_0(x) + \Tb), \quad \forall
	x\in\Omega\bs\Sigma.
\end{equation}
Now, let $\phi\compose\vp'=(u',v')$, and let
\[
	T = \max\{ \, \sup_{\D\Omega} u + 1, \, \Tb, \, (2a)^{-1}\log2\}.
\]
Note that $\sup_{\D\Omega} u=\sup_{\D\Omega} f_{-\gamma}\compose\psi$, hence
$T$ depends only on $\psi$ and $\vp_0$.
Define
\[
	u'' = u_0 + \min\{\, u'-u_0, \, T\},
\]
and define $\vp''$ by $\phi\compose\vp''=(u'',v')$.  Then $\vp''=\vp'$ on
$\D\Omega$, hence $\vp''\in\H$.  Also, from the fact that $Q_{\vp'}(\nabla
v')$ is non-decreasing in $u'$, we have $F(\vp'') \leq
F(\vp') \leq F(\vp)$.  Furthermore, the estimate:
\[
	u'' - u_0 \leq T,
\]
holds throughout $\Omega$.  This estimate can be rewritten as
$f_{-\gamma}\compose\vp''\leq u_0 + T$, or equivalently as:
\begin{equation}	\label{u:estimate}
	\vp''(x) \in \B_{-\gamma}(u_0(x) + T), \quad \forall
	x\in \Omega\bs\Sigma.
\end{equation}
We claim that:
\begin{equation}	\label{new:ub:estimate}
	\vp''(x) \in \B_\gamma(-u_0(x) + T), \quad \forall
	x\in\Omega\bs\Sigma
\end{equation}
holds still.
Indeed, if $x\in\Omega\bs\Sigma$ is such that $u'(x)-u_0(x)\leq T$,
then $\vp''(x) = \vp'(x)$, hence~\eqref{new:ub:estimate} follows directly
from~\eqref{ub:estimate}.  On the other hand, if
$x\in\Omega\bs\Sigma$ is such that $u'(x) - u_0(x) > T$.
Then, taking
$t_0=u_0(x)$, we find:
\[
	\vp'(x) \in \B_\gamma(-t_0+T)\bs\B_{-\gamma}(t_0+T).
\]
Let $\beta=-\gamma_{v'(x)}$ be the geodesic through $\vp'(x)$ to
$\gamma(-\infty)$, then $\vp''(x)$ is the
unique point where $\beta$ intersects
$\S_{-\gamma}(t_0+T)$.  Thus, to prove the claim, it suffices to check that
$\beta$ enters $\B_{-\gamma}(t_0+T)$ before it leaves $\B_\gamma(-t_0+T)$.
The claim now follows easily from Lemma~\ref{monotone}.  Indeed,
suppose to the contrary that $\beta$ left $\B_\gamma(-t_0+T)$ before it
entered $\B_{-\gamma}(t_0+T)$, then it would do so at a point
$p\in\S_\gamma(-t_0+T)\bs\B_{-\gamma}(t_0+T)$.  However, at such a point
$p$, it follows from Lemma~\ref{monotone} that $\dot\beta$ points inwards
into $\B_\gamma(-t_0+T)$, a contradiction.
Now, from~\eqref{u:estimate} and~\eqref{new:ub:estimate},
and Lemma~\ref{bound}, we deduce that
\[
	\dist(\vp''(x),\vp_0(x)) \leq R, \quad \forall x\in\Omega\bs\Sigma,
\]
where $R=T + 2^{-1}\log 2$.  We conclude that $\vp''\in\Hc$.
This completes the proof of Lemma~\ref{lemma:inf} and of
Proposition~\ref{prop:single}.
\end{pf*}

\noindent{\em Remark\/}.  We note here that the a priori estimate
$\dist(\vp,\vp_0)\leq R$ implies that $v\to0$ on $\Sigma$.
An interesting analytic
question is whether the function $u-u_0$ is continuous and
perhaps yet smoother in a neighborhood of $\Sigma$.  The interest stems from
the fact that the ellipticity of the equations degenerates near $\Sigma$.
When $\M=\HR$, this regularity question
was studied in the axially symmetric case
in~\cite{weinstein90,weinstein92}, and in the general
case in~\cite{li92,li93}.  This question is not addressed here.

\section{The Case $N\geq2$} \label{multiple}

In this section, we finish the proof of Theorem~\ref{main} by proving:

\begin{prop}	\label{prop:multiple}
Let $N\geq2$, and for $1\leq i\leq N$
let $\vp_i$ be a $\Sigma_i$-singular map
into $\gamma_i$, and let $\psi\in C^{2,\alpha}(\D\Omega;M)$.
Then, there exists a unique
harmonic map $\vp\in C^\infty(\Omega\bs\Sigma;M)\cap
C^{2,\alpha}(\overline\Omega\bs\Sigma;M)$ such that
$\vp=\psi$ on $\D\Omega$, and $\vp$ is asymptotic to $\vp_i$ near
$\Sigma_i$ for each $1\leq i\leq N$.
\end{prop}
\noindent{\em Proof\/}.
The proof of Proposition~\ref{prop:multiple} follows the same
outline as the proof of Proposition~\ref{prop:single},
but there are a few more technical points.
Let $\{\Omega_i\}_{i=1}^N$ be an open cover of $\Omega$ such that
$\Sigma_i\subset\Omega_i$, and $\Sigma_i\cap\Omega_{i'}=\emptyset$
for $i\ne i'$.  Suppose $\vp$ and $\vp'$ are harmonic maps which are
asymptotic to $\vp_i$ near $\Sigma_i$ for each $1\leq i\leq N$,
and agree with $\psi$ on $\D\Omega$.
Then, $\dist(\vp,\vp')\in L^\infty(\Omega_i\bs\Sigma_i)$
for each $1\leq i\leq N$, hence $\dist(\vp,\vp')\in
L^\infty(\Omega\bs\Sigma)$.  Thus, as in the proof of
Proposition~\ref{prop:single}, it follows that $\vp=\vp'$.

We now turn to the proof of existence.  We may, without loss of generality,
assume that all the geodesics
$\gamma_i$ have the same initial point $\gamma_1(-\infty)\in\D M$, and are
parameterized so that $f_{-\gamma_i}=f_{-\gamma_1}$ for
all $1\leq i\leq N$.
Let $u=f_{-\gamma_1}$, and let $\phi=(u,v)$ be the corresponding coordinate
system given in Lemma~\ref{lemma:Qequiv}.
We will as before identify $\vp$ and its
parameterization $\phi\compose\vp=(u,v)$ where no confusion arises.
Let $\vp_i=\gamma_i\compose u_i$, then
$\vp_i=(u_i,w_i)$, where $w_i\in\R$ are constants,
$\Delta u_i=0$ on $\Omega\bs\Sigma_i$,
and we assume without loss of generality that $u_i=0$ on $\D\Omega$, and
consequently $u_i>0$ in $\Omega\bs\Sigma_i$.  In addition, we know that there
exists a constant $\delta>0$ such that
\[
	\norm{\nabla u_i}^2 = \norm{d\vp_i}^2 \geq \delta r_i^{-2},
\]
in a neighborhood of $\Sigma_i$,
where $r_i(x) = \dist(x,\Sigma_i)$.  Set $u_0=\sum_{i=1}^N u_i$.

Extend $\psi$ to a map $\tilde\psi\in C^\infty(\Omega;M)\cap
C^{2,\alpha}(\overline\Omega;M)$ which for each $1\leq i\leq N$
maps a neighborhood of $\Sigma_i$ to the point $\gamma_i(0)$, and write
$(\tilde u,\tilde v)=\tilde\psi$.  Then $(\tilde u,\tilde v)=(0,w_i)$ in a
neighborhood of $\Sigma_i$.  Define $\H$ to be the space of maps
$\vp=(u,v)\colon\Omega\bs\Sigma\to M$ such that $u-u_0-\tilde u\in
H_{1,0}(\Omega)$, $v-\tilde v\in \cap_{i=1}^N
H_{1,0}^{\vp_i}(\Omega;\R^{m-1})$, and $\dist(\vp,\vp_i)\in
L^\infty(\Omega_i\bs\Sigma_i)$ for each $1\leq i\leq N$.
For maps $\vp\in\H$, and
$\Omega'\subset\Omega$ define $F_{\Omega'}$ and $F$ as in
Proposition~\ref{prop:single}:
\[
	F_{\Omega'}(\vp) = \int_{\Omega} \left\{
	\norm{\nabla(u-u_0)}^2 + Q_\vp(\nabla v)\right\},
\]
and $F=F_\Omega$.
Then, as before, if $\vp\in\H$ is a
minimizer of $F$, then $\vp$ is a harmonic map on $\Omega\bs\Sigma$, hence
$\vp\in C^\infty(\Omega\bs\Sigma;M)\cap
C^{2,\alpha}(\overline\Omega\bs\Sigma;M)$, and $\vp$ is asymptotic to
$\vp_i$ near $\Sigma_i$ for each $1\leq i\leq N$.  Thus, to prove
Proposition~\ref{prop:multiple}, it suffices to show that $F$ has a
minimizer $\vp\in\H$.

For any $R>0$ define the space $\H_R$ of maps $\vp\in\H$ which for each
$1\leq i \leq N$ satisfy $\dist(\vp,\vp_i)\leq R$
for a.e.\ $x\in\Omega_i\bs\Sigma_i$.
We first show that $F$ has a minimizer in $\H_R$.
For this purpose, we need the following lemma:

\begin{lemma}	\label{lemma:Nequiv}
Let $R>0$, then there is $c>0$ such that for all $\vp\in\Hc$, there holds:
\begin{equation}	\label{Nequivalent}
	\frac1c\, Q_{\vp_i(x)}(\xi) \leq Q_{\vp(x)}(\xi) \leq c\,
	Q_{\vp_i(x)}(\xi), \quad \forall\xi\in\R^{m-1}, \> \text{a.e.\
	$x\in\Omega_i\bs\Sigma_i$}.
\end{equation}
\end{lemma}
\begin{pf*}{Proof of Lemma~\ref{lemma:Nequiv}}
The lemma follows immediately from (i) in Lemma~\ref{lemma:Qequiv},
see the proof of Lemma~\ref{lemma:equiv}.
\end{pf*}

Now, let $\vp'_j=(u'_j,v'_j)\in\Hc$ be a minimizing sequence.  Then, as
before, $u'_j-u_0$ is bounded in $H_1(\Omega)$, and we may assume that it
converges weakly and pointwise a.e.\ in $\Omega$ to $u-u_0\in H_1(\Omega)$.
Thus,
\[
	\int_\Omega \norm{\nabla(u-u_0)}^2 \leq \liminf \int_\Omega
	\norm{\nabla(u'_j-u_0)}^2,
\]
and $u-u_0-\tilde u\in H_{1,0}(\Omega)$.  We claim that $v'_j-\tilde v$
is bounded in $H_{1,0}^{\vp_i}(\Omega;\R^{m-1})$ for each $1\leq i\leq N$.
To see this, note first that in view of Lemma~\ref{lemma:Nequiv}, we have
\begin{align*}
	\int_{\Omega_i} Q_{\vp_i}(\nabla v'_j - \nabla\tilde v)
	& \leq
	2 c\, \int_{\Omega_i} Q_{\vp'_j}(\nabla v'_j) +
	2 \int_{\Omega_i} Q_{\vp_i}(\nabla\tilde v)\\
	&\leq 2c\, (\, \inf_{\Hc} F + c_2) + 2 \int_{\Omega_i}
	Q_{\vp_i}(\nabla\tilde v)\\
	&\leq 2c\, (\inf_{\Hc} F + c_2) + 2 \max_{1\leq i\leq N}
	\int_{\Omega} Q_{\vp_i} (\nabla \tilde v),
\end{align*}
for some $c_2\geq0$.  Denote the constant on the right-hand side by $c_3$.
Now $u_i$ is bounded in $\Omega_{i'}\bs\Sigma_{i'}$
for $i'\ne i$.  Let $t_0 =
\max\{ \sup_{\Omega_{i'}\bs\Sigma_{i'}} u_i: \> i\ne i'\}$,
then for any $x\in\Omega_{i'}\bs\Sigma_{i'}$,
and any
$\xi\in\R^{m-1}$, we have:
\begin{alignat*}{2}
	&Q_{\vp_i(x)}(\xi) &\leq Q_{\gamma_i(t_0)}(\xi) &\leq c_4 \,
	\norm{\xi}^2 \\
	&Q_{\vp_{i'}(x)}(\xi) \,
	&\geq Q_{\gamma_{i'}(0)}(\xi) &\geq c_4^{-1} \,
	\norm{\xi}^2,
\end{alignat*}
for some $c_4\geq1$.
Therefore, we obtain:
\[
	\int_{\Omega_{i'}} Q_{\vp_i}(\nabla v'_j - \nabla \tilde v) \leq
	c_4^2
	\int_{\Omega_{i'}} Q_{\vp_{i'}}(\nabla v'_j -
	\nabla\tilde v) \leq c_4^2 \, c_3,
\]
and it follows that:
\[
	\int_{\Omega} Q_{\vp_i}(\nabla v'_j - \nabla\tilde v) \leq
	\sum_{i=1}^N \int_{\Omega_i}
	Q_{\vp_i}(\nabla v'_j - \nabla\tilde v) \leq
	c_3 + (N-1) \,c_3 \,c_4^2.
\]
Since the norm on the left hand side is equivalent to the full norm on
$H_{1,0}^{\vp_i}(\Omega;\R^{m-1})$, we obtain that indeed $v'_j-\tilde v$ is
bounded in this space.  Thus, we may assume that $v'_j-\tilde v$ converges
weakly in $H_{1,0}^{\vp_i}(\Omega;\R^{m-1})$ for each $1\leq i\leq N$,
and pointwise a.e.\ in $\Omega$ to $v-\tilde v\in
H_{1,0}^{\vp_i}(\Omega;\R^{m-1})$.  It follows at
once that $\vp\in\H$, and furthermore, since $\vp'_j(x)\to\vp(x)$ for a.e.\
$x\in\Omega$, we obtain that $\vp\in\Hc$.  We now claim that
\[
	\int_\Omega Q_\vp(\nabla v) = \lim \int_\Omega Q_\vp(\nabla v,
	\nabla v'_j).
\]
Indeed, let $\{\chi_i\}_{i=1}^N$
be a partition of unity subordinate to the cover
$\{\Omega_i\}$.  Then, using Lemma~\ref{lemma:Nequiv}, we obtain:
\begin{align*}
	\norm{ \int_\Omega Q_\vp(\nabla v) - \int_\Omega Q_\vp(\nabla
	v,\nabla v'_j)} &=
	\norm{ \int_\Omega \sum_{i=1}^N
	Q_\vp\bigl(\nabla(\chi_i v), \nabla v - \nabla
	v'_j\bigr)} \\
	& \leq \sum_{i=1}^N \norm{ \int_\Omega Q_\vp\bigl(\nabla (\chi_i v),
	\nabla v - \nabla v'_j \bigr)} \\
	&\leq c \, \sum_{i=1}^N
	\norm{\int_{\Omega_i} Q_{\vp_i}\bigl(\nabla (\chi_i v), \nabla v -
	\nabla v'_j\bigr)},
\end{align*}
which tends to zero since certainly $v-v'_j$ tends to zero weakly in
$H_{1,0}^{\vp_i}(\Omega_i;\R^{m-1})$.  Now, the argument in the proof of
Proposition~\ref{prop:single} applies,
and we deduce that $\vp$ is a minimizer of $F$
on $\Hc$.

Once more, it remains only to prove that some $R>0$ large enough
$\inf_{\H} F = \inf_{\Hc} F$.
As before, we state this as a lemma.

\begin{lemma}	\label{lemma:Ninf}
There is a constant $R>0$ such that for every $\epsilon>0$, and every
$\vp\in\H$, there is $\vp'\in\Hc$ such that $F(\vp') \leq F(\vp) +
\epsilon$.
\end{lemma}
\begin{pf*}{Proof of Lemma~\ref{lemma:Ninf}}
Let $\H^*$ be the space of maps $\vp\in\H$ such that
$v=w_i$ in a neighborhood
of $\Sigma_i$ for each $i$.
The proof of the lemma is divided into steps
(i) and (ii) as in the proof of Lemma~\ref{lemma:inf}.  The proof of (i) is
practically unchanged.  We immediately turn to the proof of (ii).
Introduce the new coordinate system $\phib=(\ub,\vb)$, where
$\ub=f_{\gamma_1}$.  Again, write the metric on $\M$ as
\[
	ds^2 = d\ub\,{}^2 + \Qb_p(d\vb).
\]
Set $\ub_0 = - u_1 + \sum_{i=2}^N u_i$.
Let $\vp\in\H^*$, and write $\phib\compose\vp=(\ub,\vb)$.
Then in $\Omega\bs\Sigma$, we find:
\begin{align*}
	\norm{\nabla(u-u_0)}^2 + Q_\vp(\nabla v)
	&= \norm{\nabla(\ub - \ub_0)}^2 + \Qb_\vp(\nabla\vb) \\
	  &\qquad\qquad\qquad
	     - 2\,\nabla u \cdot \nabla u_0
	     + 2\,\nabla \ub \cdot \nabla \ub_0
	     + \norm{\nabla u_0}^2 - \norm{\nabla \ub_0}^2\\
	&= \norm{\nabla(\ub - \ub_0)}^2 + \Qb_\vp(\nabla\vb)
	 - 2\,\nabla(u+\ub) \cdot \nabla u_1 \\
	 & \qquad\qquad\qquad
	     - 2 \textstyle\sum_{i=2}^N \nabla (u-\ub) \cdot \nabla u_i
	     + 4 \textstyle\sum_{i=2}^N \nabla u_1 \cdot \nabla u_i.
\end{align*}
We wish to integrate this identity over $\Omega$ in order to obtain an
integral identity analogous to~\eqref{integral:identity}.
Some care must be taken due to the singularities at $\Sigma_i$.
Integrate first over $\Omega^\epsilon=
\{x\in\Omega:\dist(x,\Sigma)>\epsilon\}$:
\begin{equation}	\label{identity:N:epsilon}
	F_{\Omega^\epsilon}(\vp) = \int_{\Omega^\epsilon}
	\left\{ \norm{\nabla(\ub - \ub_0)}^2 + \Qb_\vp(\nabla\vb) \right\}
	+ I_1(\epsilon) + I_2(\epsilon) + I_3(\epsilon),
\end{equation}
where we have set:
\begin{align}
	\label{i1}
	I_1(\epsilon) &= -2 \int_{\Omega^\epsilon} \nabla(u+\ub) \cdot
	\nabla u_1 =
	-2 \int_{\Omega^\epsilon} \Div\bigl( (u+\ub)\, \nabla
	u_1\bigr) \\
	\label{i2}
	I_2(\epsilon) &= -2\sum_{i=2}^N \int_{\Omega^\epsilon}
	\nabla(u-\ub)\cdot \nabla u_i =
	-2 \sum_{i=2}^N
	\int_{\Omega^\epsilon} \Div\bigl( (u-\ub) \,\nabla
	u_i\bigr) \\
	\label{i3}
	I_3(\epsilon) &= 4 \sum_{i=2}^N \int_{\Omega^\epsilon}
	\nabla u_1 \cdot \nabla u_i.
\end{align}
For $\epsilon>0$
small enough, we can decompose the boundary $\D\Omega^\epsilon$ into
a disjoint union $\cup_{i=1}^N \D_i^\epsilon \cup\D\Omega$, where
$\D_i^\epsilon=\{x\in\Omega: \dist(x,\Sigma_i)=\epsilon\}\subset\Omega_i$.
Integrate~\eqref{i1} by parts:
\[
	I_1(\epsilon) = -2 \int_{\D\Omega} (u +\ub)\, \frac{\D
	u_1}{\D n} - 2 \int_{\D^\epsilon_1} (u + \ub)\, \frac{\D
	u_1}{\D n} -
	2 \sum_{i=2}^N \int_{\D_i^\epsilon} (u + \ub)\, \frac{\D
	u_1}{\D n}.
\]
For $\epsilon$ small
enough, $v=0$, and hence $u+\ub=0$ on $\D_1^\epsilon$,
so the second term above vanishes.  Let $i\geq2$,
and write $\ub_i=f_{\gamma_1}\compose\vp_i$.
Then, since $\vp\in\H_{R'}$ for some $R'>0$, we have almost
everywhere in $\Omega_i\bs\Sigma_i$:
\begin{align*}
	\norm{u-u_i} &\leq R'\\
	\norm{\ub-\ub_i} &\leq R'.
\end{align*}
{}From Lemma~\ref{constant}, Equation~\eqref{minus}, we
have $u_i - \ub_i\to d_i$ as $x\to\Sigma_i$, for some $d_i\in\R$, and
from~\eqref{gamma:i}, we have $\norm{u_i-\ub_i}\leq D$, where
\[
	D = \max_{2\leq i\leq N}
	\bigl\{\norm{d_i}, \norm{f_{-\gamma_1}\compose\gamma_i(0) -
	f_{\gamma_1}\compose\gamma_i(0)}\bigr\}.
\]
Furthermore, there is a constant $c_6>0$ such that
$\norm{\D u_1/\D n}\leq c_6$
in $\Omega_i\bs\Sigma_i$.  Consequently, taking $a_i(\epsilon)$ to be the
$(n-1)$-Hausdorff measure of $\D^\epsilon_i$,
we can estimate:
\begin{equation}	\label{boundary:integral}
	\begin{split}
	\norm{\int_{\D_i^\epsilon} (u + \ub)\, \frac{\D u_1}{\D n}}
	& \leq c_6
	\int_{\D_i^\epsilon} \biggl\{ \norm{u_i-\ub_i} + 2u_i + 2R'
	\biggr\} \\
	& \leq
	\biggl( (2R'+D) \, a_i(\epsilon) +
	2 \int_{\D^\epsilon_i} u_i \biggr)\to 0,
	\end{split}
\end{equation}
as $\epsilon\to0$, by Lemma~\ref{lemma:integral:estimate}.
We deduce that, as $\epsilon\to0$:
\begin{equation}	\label{k1}
	I_1(\epsilon) \to K_1 = -2 \int_{\D\Omega} (u+\ub)\,
	\frac{\D u_1}{\D n},
\end{equation}
which clearly depends only on $\psi$ and $\vp_1$.
Now, we have
\[
	I_2(\epsilon) = -2\sum_{i=2}^N \int_{\D\Omega} (u-\ub) \,
	\frac{\D u_i}{\D n} - 2\sum_{i=2}^N \sum_{i'=1}^N
	\int_{\D^\epsilon_{i'}} (u-\ub)\, \frac{\D u_i}{\D n}.
\]
An estimate similar to~\eqref{boundary:integral} shows that all the terms in
the second sum are zero except when $i'=i\geq2$.  To compute this term,
let:
\[
	e_i = \int_{\D^\epsilon_i} \frac{\D u_i}{\D n},
\]
be the charge of $u_i$, see Lemma~\ref{lemma:integral:estimate}.
If $\epsilon>0$ is small enough,
we have for each $x\in\D^\epsilon_i$
that $\vp(x)$ is a point along $\gamma_i$, and $u(x)\geq u_i(x)-R'$ for
a.e.\ $x\in\Omega_i\bs\Sigma_i$.  It
follows from Lemma~\ref{constant}, Equation~\eqref{minus}, that
$(u-\ub)|\D^\epsilon_i\to d_i$ as $\epsilon\to0$, and $\norm{u-\ub}
\leq \norm{u_i-\ub_i} + 2R' \leq D + 2R'$.  Also, for $\epsilon>0$ small
enough $(\D u_i/\D n) >0$, and clearly,
\[
	\int_{\D^\epsilon_i} (D+2R') \, \frac{\D u_i}{\D n} = (D +2R')\, e_i,
\]
for each $\epsilon>0$. Thus, integrating over a representative surface
$\D^{\epsilon_0}_i$, and using Fatou's Lemma, we conclude that
\begin{equation}	\label{di}
	\int_{\D^\epsilon_i} (u - \ub)\, \frac{\D u_i}{\D n} \to
	d_i \, e_i,
\end{equation}
as $\epsilon\to0$.  Summing up, we have obtained that, as $\epsilon\to0$:
\begin{equation}	\label{k2}
	I_2(\epsilon) \to K_2 = -2 \sum_{i=2}^N
	\left\{ \int_{\D\Omega} (u-\ub)\, \frac{\D u_i}{\D n} +
	d_i \, e_i \right\},
\end{equation}
which clearly depends only on $\psi$ and $\vp_i$, $2\leq i\leq N$.
Now, for $\epsilon>0$ small enough,
and $i\geq2$, we calculate:
\begin{align*}
	\int_{\Omega^\epsilon} \nabla u_1 \cdot \nabla u_i
	&= \int_{\Omega^\epsilon\cap\Omega_1}
	\Div( u_i\, \nabla u_1) +
	\int_{\Omega^\epsilon\bs\Omega_1}
	\Div( u_1\, \nabla u_i) \\
	&= \int_{\D_1} u_i\, \frac{\D u_1}{\D n}
	+ \int_{\D^\epsilon_1} u_i \frac{\D u_1}{\D n} +
	\int_{\D'_1} u_1\, \frac{\D u_i}{\D n}
	+ \int_{\D^\epsilon_i} u_1 \frac{\D u_i}{\D n}
	+ \sum_{1\ne i'\ne i} \int_{\D^\epsilon_{i'}} u_1\,
	\frac{\D u_i}{\D n},
\end{align*}
where $\D_1=\D(\Omega\cap\Omega_1)$, and $\D'_1=\D(\Omega\bs\Omega_1)$.
Since $u_1$ and $u_i$ are smooth in $\Omega_{i'}$ for $1\ne i'\ne i$, all
the terms in the last sum tend to zero as $\epsilon\to0$.
In $\Omega_1$, $u_i$ is smooth, thus by an argument analogous to the one
leading to~\eqref{di}, we have that, as $\epsilon\to0$:
\[
	\int_{\D^\epsilon_1} u_i \, \frac{\D u_1}{\D n} \to k_1,
\]
for some $k_1\in\R$.  Similarly:
\[
	\int_{\D^\epsilon_i} u_1 \, \frac{\D u_i}{\D n} \to k_i,
\]
for some $k_i\in\R$.  Summing up, we have obtained that:
\begin{equation}	\label{k3}
	I_3(\epsilon) \to K_3 = k_1 + \sum_{i=2}^N \left(
	\int_{\D_1} u_i\, \frac{\D u_1}{\D n} +
	\int_{\D'_1} u_1\, \frac{\D u_i}{\D n} +
	k_i \right).
\end{equation}
Combining~\eqref{identity:N:epsilon} with~\eqref{k1}, \eqref{k2}
and~\eqref{k3}, we conclude that
\begin{equation}	\label{integral:N:identity}
	F(\vp) = \int_{\Omega^\epsilon}
	\left\{ \norm{\nabla(\ub - \ub_0)}^2 + \Qb_\vp(\nabla\vb) \right\}
	+ K,
\end{equation}
where $K=K_1+K_2+K_3$ is independent of $\vp\in\H^*$.  Now the proof can
proceed as in Lemma~\ref{lemma:inf}.  Truncate $\ub-\ub_0$ above at
\[
	\Tb_1 = \sup_{\D\Omega} \ub + 1,
\]
to get a map $\vp'=(u',v')\in\H$ satisfying $F(\vp')\leq F(\vp)$, and
\begin{equation}	\label{ub:N:estimate}
	\vp'(x) \in \B_{\gamma_1}(\ub_0(x) + \Tb_1), \quad \forall
	x\in\Omega\bs\Sigma.
\end{equation}
Let $c_7=\sup_{\Omega_1} \sum_{i=2}^N u_i$, then from~\eqref{ub:N:estimate},
we obtain:
\begin{equation}	\label{ub:N:estimate:1}
	\vp'(x) \in \B_{\gamma_1}(-u_1(x) + c_7 + \Tb_1), \quad \forall
	x\in\Omega_1\bs\Sigma_1.
\end{equation}
Now, truncate $u'-u_0$ above at
\[
	T_1 = \max\{ \sup_{\D\Omega} u+1, \Tb_1, (2a)^{-1}\log2\},
\]
to get a map $\vp''=(u'',v')\in\H$ satisfying $F(\vp'')\leq F(\vp')\leq
F(\vp)$, and
\begin{equation}	\label{u:N:estimate}
	\vp''(x)\in\B_{-\gamma_1}(u_1(x)+T'_1),
	\quad \forall x\in\Omega\bs\Sigma,
\end{equation}
where $T'_1=T_1+c_7$.
As in the proof of Lemma~\ref{lemma:inf}, it follows from
Lemma~\ref{monotone} and~\eqref{ub:N:estimate:1} that:
\begin{equation}	\label{ub:N:new:estimate:1}
	\vp''(x) \in \B_{\gamma_1}(-u_1(x) + T'_1), \quad \forall
	x\in\Omega_1\bs\Sigma_1.
\end{equation}
Therefore, combining~\eqref{u:N:estimate} with ~\eqref{ub:N:new:estimate:1},
and Lemma~\ref{bound}, we conclude that:
\begin{equation}	\label{vp:N:estimate}
	\vp''(x) \in B_{R_1}\bigl(\vp_1(x)\bigr), \quad \forall
	x\in\Omega_1\bs\Sigma_1,
\end{equation}
where $R_1 = T'_1 + a^{-1}\log2$.
Now consider the map $\vp''|\Omega'$, where $\Omega'=\cup_{i=2}^N\Omega_i$
and note that, since $\D\Omega'=(\D\Omega\cap\D\Omega') \cup
(\D\Omega'\cap\Omega_1)$,
\eqref{vp:N:estimate} together with $\psi$
give a pointwise a priori estimate for $\vp''$
throughout $\D\Omega'$.  Thus one can proceed by induction to
obtain a map $\vp'''\in\H$ which satisfies $F(\vp''')\leq F(\vp)$, and
for each $1\leq i\leq N$:
\[
	\vp'''(x) \in B_{R_i} \bigl(\vp_i(x)\bigr), \quad \forall
	x\in\Omega_i\bs\Sigma_i,
\]
for some
constants $R_i$ depending only on the boundary map $\psi$, and
the $\Sigma_i$-singular harmonic maps $\vp_i$.
Set $R=\max_i R_i$, then we have obtained $\vp'''\in\Hc$ with $F(\vp''')\leq
F(\vp)$.
This completes the proof of Lemma~\ref{lemma:Ninf}, and of
Proposition~\ref{prop:multiple}.
\end{pf*}
\vspace{1em}

\pagebreak

\appendix
\begin{center}
	\sc Appendix
\end{center}
\hspace{\parindent} In this appendix, we derive
Equations~\eqref{realhyperbolic}--\eqref{quaternionhyperbolic}
and~\eqref{busemann} used in
Lemma~\ref{lemma:Qequiv}.
Our starting point will be the disk models for the spaces
$\Ha^\ell_{\K}$, as given in~\cite[Section 19]{mostow}, but first
we briefly recall how these models are obtained.
Consider $\K^{\ell+1}$ as a right $\K$-module equipped with the bilinear
form:
\begin{equation}	\label{bilinear}
	\inner{\bx}{\by} = \xb_0 y_0 - \sum_{k=1}^{\ell} \xb_k y_k.
\end{equation}
where $x \mapsto \xb$ is the standard involution of $\K$.
Let $\pi\colon\K^{\ell+1}\to \P\K^\ell$ denote the natural map onto the
projective space $\P\K^{\ell}$ of $\K^{\ell+1}$.
Consider $\overline{\Do} = \{ \bx\in\K^{\ell+1}: \inner{\bx}{\bx} >
0 \}$, and its image $\Do = \pi(\overline{\Do})$ under $\pi$.  Introduce
non-homogeneous coordinates $\bz\in\K^\ell$ on $\Do$ by
\begin{equation}	\label{non-homogeneous}
	z_k = x_k x_0^{-1}, \qquad 1\leq k\leq \ell,
\end{equation}
and define the bilinear form:
\[
	\scalar{\bz}{\bw} = \sum_{k=1}^\ell \zb_k w_k.
\]
It is easily seen that
$\Do=\{\bz\in\K^\ell: \scalar{\bz}{\bz} < 1 \}$.
Now, the metric
\begin{equation}	\label{ds2non-homogeneous}
	d \overline{s}\,{}^2 = - \frac{\inner{d\bx}{d\bx}}{\inner{\bx}{\bx}}
	+ \frac{\norm{\inner{d\bx}{\bx}}^2}{\inner{\bx}{\bx}^{2}}
\end{equation}
on $\overline \Do$
is invariant under the right action of $\K$, hence it
induces a quotient metric on $\Do$, which
in terms of the coordinates $\bz$ turns out to be:
\begin{equation}	\label{ds2}
	ds^2 = \frac{\norm{d\bz}^2}{1 - \norm{\bz}^2} +
		\frac{\norm{\scalar{d\bz}{\bz}}^2}{
		\bigl(1-\norm{\bz}^2\bigr)^2}.
\end{equation}
To see this, use~\eqref{non-homogeneous} to get $dz_k = dx_k x_0^{-1} -
x_k x_0^{-1} dx_0 x_0^{-1}$, and substitute into the left-hand side
of~\eqref{ds2} to obtain the left-hand side of~\eqref{ds2non-homogeneous}.
If we now let $GL(\ell,\K)$ act on the left on $\K^{\ell+1}$,
and we let $G\subset GL(\ell,\K)$ be a subgroup which leaves the
form~\eqref{bilinear} invariant, then
clearly, $G$ leaves $\overline{\Do}$ invariant.
Also, the action of $G$ commutes with the right action of $\K$ on
$\K^{\ell+1}$, hence $G$ acts on $\P\K^{\ell}$ leaving $\Do$
invariant.  Thus, we obtain a left $G$-action on $\Do$,
which leaves the metric~\eqref{ds2}
invariant.  Let $h=(h_{kj})\in G$, then this action
$(h,\bz)\mapsto h\cdot\bz$ is easily seen to be
given by the M\"obius transformation:
\begin{equation}	\label{action}
	(h\cdot\bz)_k =
	\biggl( h_{k0} + \sum_{j=1}^{\ell} h_{kj} z_j \biggr)
	\biggl( h_{00} + \sum_{j=1}^{\ell} h_{0j} z_j \biggr)^{-1}.
\end{equation}
Let $G=SO_0(1,\ell)$ when $\K=\R$,
$G=SU(1,\ell)$ when $\K=\Co$, and $G=Sp(1,\ell)$ when $\K=\Ha$.
It is not difficult to see that $G$ is transitive on $\Do$.  Let
$\bx_0 = (1,0,\dots,0)\in\Do$, and let
$\o=\pi(\bx_0)$, then the isotropy subgroup $H_\o\subset G$
of $\o$ is readily calculated.  When $\K=\R$, we have
$H_\o=SO(\ell)$, while when
$\K=\Co$, we have $H_\o=S(U(1)\times U(\ell))$,
and when $\K=\Ha$, we have $H_\o=Sp(1)\times Sp(\ell)$.
Thus, we conclude that $(\Do,ds^2)$ can serve as a
model for $\Ha_{\K}^\ell$.  In order to calculate the distance function on
$\Ha_{\K}^\ell$, first note that the
real lines $\bz t$ through $\o$ in $\Do$ are geodesics,
hence if $\bz\in\Do$, then $\dist(\o, \bz) = \tanh^{-1}(\norm{\bz})$.
Now use the action~\eqref{action},
to obtain for each $\bz\in\Do$ an element $h\in G$ such that
$h\cdot\bz=\o$.
Since for this $h$, we have
$1 - \norm{h\cdot\bw}^2 = \norm{1-\scalar{\bz}{\bw}}^{-2}
(1-\norm{\bz}^2)(1-\norm{\bw}^2)$ for all $\bw\in\Do$, that leads to:
\begin{equation}	\label{distance}
	\dist(\bz,\bw) = \cosh^{-1} \left(
	\frac{\norm{1 -
	\scalar{\bz}{\bw}}}{\bigl(1-\norm{\bz}^2\bigr)^{1/2}
	\bigl(1-\norm{\bw}^2\bigr)^{1/2}} \right).
\end{equation}

We now wish to obtain an upper half-space model for $\Ha_{\K}^\ell$.  For
this purpose, we compute the Busemann function $f_{-\gamma}$ where
$\gamma$ is any geodesic.  Since $H_{\K}^\ell$ is
homogeneous and isotropic, we may assume that $\Do$ has been
so constructed that $\gamma$ is the real line through $\o$ tangent to
$\e_1=(1,0\dots,0)$, i.e.\ $\gamma(t)=\e_1 \! \tanh(t)$.  It is now
straightforward to check that:
\begin{equation}	\label{-busemann}
	f_{-\gamma}(\bz) = \lim_{t\to\infty}
	\bigl[ \dist\bigl(\bz,\e_1 \! \tanh(-t)\bigr) - t\bigr] =
	\log\frac{\norm{1+z_1}}{\bigl(1-\norm{\bz}^2\bigr)^{1/2}}.
\end{equation}
We set $u=f_{-\gamma}$.
Following~\cite{mazur}, we define:
\[
	w_1 = (1-z_1)(1+z_1)^{-1}, \qquad w_k = z_k (1+z_1)^{-1}, \quad
	\text{for $2\leq k\leq\ell$}.
\]
We find that:
\[
	\Re w_1 - \sum_{k=2}^\ell \norm{w_j}^2 =
	\frac{1-\norm{\bz}^2}{\norm{1+z_1}^2} =
	e^{-2u},
\]
hence:
\begin{equation}	\label{realpart}
	d(e^{-2u}) = \Re\bigl[ dw_1
	- 2\sum_{k=2}^\ell \overline w_j dw_j \bigr].
\end{equation}
Furthermore, we have:
\begin{equation}	\label{first}
\begin{split}
	\frac14 e^{4u} \bigl|
	dw_1 \! - 2 \! \sum_{k=2}^\ell \overline w_k \, dw_k \bigr|^2
	&= \frac{\norm{dz_1}^2}{\norm{1+z_1}^2} +
	\frac{\norm{\scalar{\bz}{d\bz}}^2}{\bigl(1-\norm{\bz}^2\bigr)^2}
	\\[1ex]
	&  {} + \frac{2\,
	\Re\bigl[ (1+\zb_1)^{-1} d\zb_1 (1+z_1) \scalar{\bz}{d\bz}
	(1+z_1)^{-1} \bigr]}{1-\norm{\bz}^2}.
\end{split}
\end{equation}
and
\begin{equation}
\begin{split}	\label{second}
	e^{2u} \sum_{k=2}^\ell \norm{dw_k}^2 &=
	\frac{\norm{d\bz}^2 -\norm{dz_1}^2}{1-\norm{\bz}^2}
	+ \frac{\bigl(\norm{\bz}^2 - \norm{z_1}^2\bigr) \norm{dz_1}^2}{
	\norm{1+z_1}^2\bigl(1-\norm{\bz}^2\bigr)} \\[1ex]
	&- \frac{2\, \Re\bigl[
	(1+\zb_1)^{-1} d\zb_1 (1+z_1) \bigl( \scalar{\bz}{d\bz} - \zb_1 dz_1
	\bigr) (1+z_1)^{-1} \bigr]}{1 - \norm{\bz}^2}.
\end{split}
\end{equation}
Combining~\eqref{first} and~\eqref{second}, we conclude that the
metric~\eqref{ds2} is given by:
\begin{equation}	\label{w}
	ds^2 = \frac14 e^{4u}
	\bigl| dw_1 - 2 \sum_{k=2}^\ell \overline w_k \, dw_k \bigr|^2
	+ e^{2u} \sum_{k=2}^\ell \norm{dw_k}^2.
\end{equation}
Let $d\theta$ denote the form
$dw_1-2\sum_{k=2}^{\ell}\overline w_k dw_k$, then
from Equation~\eqref{w} and~\eqref{realpart}, we now see that:
\[
	ds^2 = du^2 + \frac14 e^{4u}
	\bigl| \bigl( d\theta - \overline{d\theta} \bigr)/2 \bigr|^2
	+ e^{2u} \sum_{k=2}^\ell \norm{dw_k}^2.
\]
We can now obtain the coordinate system $\phi=(u,v)$ as claimed in
Lemma~\ref{lemma:Qequiv}.
If $\K=\R$, we set $v_k=w_k$ for $2\leq k\leq\ell$, and immediately
obtain~\eqref{realhyperbolic}.  If $\K=\Co$, we set
\[
	v_1 = \frac12 \Im w_1, \quad v_{2k} = \Re w_k, \quad
	v_{2k+1} = \Im w_k, \quad \text{for $1\leq k\leq\ell-1$},
\]
and obtain~\eqref{complexhyperbolic}.  Finally, if $\K=\Ha$, let
$\{1,\bold i,\bold j,\bold k\}$ be the standard basis of $\Ha$ over $\R$.
We set
\begin{align*}
	\bold i v_1 + \bold j v_2 + \bold k v_3 &=  \frac14
	\bigl(w_1-\overline w_1\bigr) \\
	v_{4k} + \bold i v_{4k+1} + \bold j v_{4k+2} + \bold k v_{4k+3} &=
	w_k, \quad \text{for $1\leq k\leq\ell-1$},
\end{align*}
and obtain~\eqref{quaternionhyperbolic}.

In order to derive~\eqref{busemann}, we observe that:
\[
	f_{\gamma}(\bz) = \lim_{t\to\infty} [\dist\bigl(\bz,\e_1
	\tanh(t)\bigr) - t] = \log
	\frac{\norm{1-z_1}}{(1-\norm{\bz}^2)^{1/2}},
\]
and
\[
	\norm{w_1}^{-2} \left( \Re w_1 - \sum_{k=2}^\ell \norm{w_k}^2
	\right) = \Re w^{-1} - \norm{w_1}^{-2} \sum_{k=2}^\ell \norm{w_k}^2
	= e^{-2f_\gamma}.
\]
Thus, setting $\ub=f_\gamma$, we have $e^{2\ub} = e^{2u} \norm{w_1}^2$.
Since
\[
	\norm{w_1}^2 = \left( e^{-2u} + \sum_{k=2}^\ell \norm{w_k}^2
	\right)^2 + \frac14 \norm{w_1 - \overline w_1}^2,
\]
we immediately obtain~\eqref{busemann}.

\end{document}